\documentclass[11pt]{article}

\usepackage{graphicx}
\usepackage{color,url,amssymb}
\usepackage{amsmath,mathtools,amsthm}
\usepackage{verbatim,epstopdf}
\usepackage{empheq}

\usepackage{hyperref}

\newcommand{\N}{\mathbb{N}}
\newcommand{\R}{\mathbb{R}}

\newtheorem{theorem}{Theorem}[section]
\newtheorem{proposition}[theorem]{Proposition}
\newtheorem{defn}{Definition}[section]
\newtheorem{asmp}{Assumption}[section]

\newtheorem{corol}[theorem]{Corollary}
\newtheorem{lemma}[theorem]{Lemma}
\newtheorem{conj}[theorem]{Conjecture}
\newtheorem{rmk}{Remark}[section]

\let\emptyset\varnothing

\begin{document}

\pagestyle{myheadings}

\title{
\vspace*{-20mm}
Compactification for Asymptotically Autonomous Dynamical Systems: Theory, Applications and Invariant Manifolds
}

\author{
Sebastian Wieczorek\thanks{University College Cork, School of Mathematical Sciences, Western Road, Cork, T12 XF62, Ireland}, Chun Xie,${}^*$ and  Chris K.R.T. Jones\thanks{Mathematics Department, University of North Carolina, Chapel Hill, NC 27599, USA} 
}

\date{\today}

\maketitle

\vspace*{0mm}

\section*{Abstract}
We develop a general compactification framework to facilitate analysis of
nonautonomous ODEs where nonautonomous terms decay asymptotically. 
The strategy is to compactify the problem: the phase space is augmented with a bounded 
but open dimension and then extended at one or both ends by gluing in flow-invariant 
subspaces that carry autonomous dynamics of the limit systems from infinity. 
We derive the weakest decay conditions possible for the compactified system to be 
continuously differentiable on the extended phase space.
This enables us to  use equilibria and other compact invariant sets of the limit systems from infinity to analyse 
the original nonautonomous problem in the spirit of dynamical systems theory. 
Specifically, we prove that solutions of interest are contained in unique invariant 
manifolds of saddles for the limit systems when embedded in the extended phase 
space. The uniqueness holds in the general case, that is even if the compactification 
gives rise to a centre direction and the manifolds become  centre or centre-stable manifolds.
A wide range of problems including pullback attractors, rate-induced 
critical transitions (R-tipping) and nonlinear wave solutions fit naturally 
into our framework. 

\vspace{5mm}
\noindent
Keywords: asymptotically autonomous differential equations, compactification,  nonautonomous dynamical systems, rate-induced tipping, nonlinear waves, centre manifold, centre-stable manifold.

\tableofcontents

\newpage

\section{Introduction}

Let $U$ and $V$ be open subsets $U\subseteq\R^n$ and $V\subseteq\R^d$. Consider a 
{\em nonlinear nonautonomous} differential equation
\begin{equation}
\label{eq:ode}
\frac{dx}{dt}:=\dot{x}=f(x, \Gamma(t)),
\end{equation} 
with the dependent state variable $x\in U$, independent variable $t\in \R$, $C^1$-smooth nonautonomous 
term $\Gamma:\R\to V$, and $C^1$-smooth vector field $f:U\times V\to \mathbb{R}^n$. 

For an arbitrary $\Gamma(t)$, 
the theory of nonautonomous dynamical systems is presented in \cite{KloedenRasmussen2011} which summarises work on the problem, 
discusses useful concepts such as pullback attractors, and gives some very general results on attraction 
and stability. 
The main obstacle to the analysis of nonautonomous system~(\ref{eq:ode}) is the absence of compact 
invariant sets such as equilibria, limit cycles or invariant tori. This obstacle becomes evident when the system 
is augmented with $\nu=t$ as an additional dependent variable to obtain the usual 
{\em autonomous  extended system}\footnote{
We note that in the special case $\Gamma(t)$ is generated by a known set of autonomous ODEs, 
say $\dot{\Gamma} = w(\Gamma)$, the original system can be augmented with $\Gamma$ to obtain the autonomous system
$\dot{x}= f(x,\Gamma),\,\dot{\Gamma} = w(\Gamma)$, 
defined on $U\times V\subseteq \R^{n+d}$. However, $w(\Gamma)$ is not always known and may be impractical for analysis if $d$ is large. 
}
\begin{align}
\begin{split}
\label{eq:odeext0}
 \dot{x}&= f(x,\Gamma(\nu)),\\
\dot{\nu}&= 1,
\end{split}
\end{align} 
that is defined on the phase space $U\times\R$, which is an open subset of $\R^{n+1}$. 

Our strategy is to consider asymptotically decaying $\Gamma(t)$ and compactify the system (\ref{eq:odeext0}). Compactification is not a new strategy and has been exploited in many applications cited in Sec.~\ref{sec:examples} below. However, in the context of decaying $\Gamma(t)$, it has been implemented only on a case-by-case basis for specific applications and specific types of decay~\cite{Jones1986,Alexander1990,Ashwin2012,Perryman2014,Perryman2015,Alkhayuon2018}. Our viewpoint here is to develop a general framework and derive in Theorem~\ref{th2} 
the weakest decay conditions possible for which the compactified system is sufficiently smooth to enable the construction of invariant manifolds (stable, unstable, centre, etc.) in the compactified phase space.

The invariant manifolds are key in the applications we envision in Sec.~\ref{sec:examples}.  Addressing the issue of finding optimal conditions for compactification raises the problem that the manifolds of interest may involve a centre direction and need not be pure stable or unstable manifolds. This, in turn, raises the possibility that the manifolds of interest may not be unique.  We show in Theorems~\ref{prop:saddler} and~\ref{prop:sinkr} that uniqueness holds in spite of these centre directions.

\subsection{The Basic Setting}
\label{sec:bs}

The main point of this work is that one can make further progress on  problem~\eqref{eq:ode} if $\Gamma(t)$ limits to 
a (vector) constant as $t$ tends to positive or negative infinity. Specifically,
\begin{defn}
  \label{defn:ac}
We say  $\Gamma(t)$ is {\em bi-asymptotically constant} with  future limit $\Gamma^+$
and past limit $\Gamma^-$ if 
   \begin{equation}
\label{eq:ac1}
\lim_{t\to \pm\infty} \Gamma(t)= \Gamma^\pm\in\R^d.\nonumber
\end{equation}
We say  $\Gamma(t)$ is {\em asymptotically constant} if it has a future limit but not necessarily a past limit, 
or if it has a past limit but not necessarily a future limit.
\end{defn}
The main simplification is that nonautonomous system~(\ref{eq:ode}) becomes {\em asymptotically autonomous} 
in the terminology of~\cite{Markus1956}:
\begin{align}
\begin{split}
f(x,\Gamma(t)) \to f(x,\Gamma^+)\;\;\mbox{as}\;\; t\to +\infty,\;\;\mbox{or}\;\;
f(x,\Gamma(t)) \to f(x,\Gamma^-)\;\;\mbox{as}\;\; t\to -\infty,
\end{split}\nonumber
\end{align}
and we can define the autonomous {\em future limit system}
\begin{equation}
  \label{eq:odea+}
  \dot{x}=f(x,\Gamma^+),
\end{equation} 
or the autonomous {\em past limit system}
\begin{equation}
  \label{eq:odea-}
  \dot{x}=f(x,\Gamma^-).
\end{equation} 
The autonomous dynamics of the future~\eqref{eq:odea+} or past~\eqref{eq:odea-} 
limit systems will typically include equilibria and other compact invariant sets. 
Note, however, that the flow of~\eqref{eq:ode} does not contain the autonomous dynamics of~\eqref{eq:odea+} 
or~\eqref{eq:odea-} because they only appear as $t$ tends to positive and negative infinity.

The key idea is to overcome the main obstacle by using the autonomous 
dynamics and compact invariant sets of the limit systems to analyse the nonautonomous 
system~\eqref{eq:ode}. 
Specifically, we augment the nonautonomous system~\eqref{eq:ode} with an additional 
dimension that is bounded but open, and then glue to its open ends the future limit system~\eqref{eq:odea+} from $+\infty$, the past limit system~\eqref{eq:odea-} from $-\infty$, or both. In this way, we combine systems~\eqref{eq:ode},~\eqref{eq:odea+} and~\eqref{eq:odea-} into one {\em compactified system} that
is autonomous, but not necessarily  $C^1$-smooth.
The compactification approach presented in this paper is  
akin to the Poincar\'e-type  phase-space compactification  
(of the $x$-dimensions) that enables analysis of 
dynamical behaviour at infinity~\cite{Delgado1995,Krauskopf1997,Dumortier2005,Messias2009,Matsue2018,Giraldo2018},
collisions in many-body problems~\cite{McGehee1974}, stability of nonlinear waves~\cite{Jones1986,Alexander1990}, as well as higher-codimension bifurcations~\cite{Dumortier1993} and 
canard solutions in slow-fast systems~\cite{Krupa2001,Szmolyan2001,Kosiuk2011,Kuehn2014,Kristiansen2017} 
via blow up of singularities of vector fields. The difference between our work and these studies
is twofold. Firstly, we compactify the augmented $\nu$-dimension, but not necessarily the $x$-dimensions.
Secondly, we consider a  general case with arbitrary decay of $\Gamma(t)$ and
\begin{itemize}
    \item
    Derive in Theorem~\ref{th2} the optimal conditions for the compatification: the {\em weakest decay of $\dot{\Gamma}(t)$ possible} that allows us to construct a {\em $C^1$-smooth compactified system}. 
    
    \item 
    Prove in Theorems~\ref{prop:saddler} and~\ref{prop:sinkr} that the geometric object of interest, namely the set of orbits tending to a saddle at infinity, is a {\em unique invariant  manifold} in the compactified phase space, even if a  centre direction arises from non-exponential decay of $\Gamma(t)$.
\end{itemize}
Our framework enables analysis of asymptotically autonomous dynamical systems~\cite{Markus1956} in terms of unique invariant manifolds of 
saddles from infinity in a general setting, where the decay of $\Gamma(t)$ ranges from super-exponential to sub-logarithmic. Thus, it provides an alternative tool and complements the existing approaches  based on pullback attractors~\cite{Langa2003,Rasmussen2008} and asymptotic equivalence of the nonautonomous system and the autonomous limit systems~\cite{Thieme1994,Robinson1996}.

The paper is organised as follows.
In Sec.~\ref{sec:examples} we describe different problems from applications that provide direct 
motivation for this study and fit naturally into the framework of~\eqref{eq:ode} with asymptotically constant $\Gamma(t)$.
In Sec.~\ref{sec:compact} we develop the general compactification technique and introduce the 
underpinning coordinate transformation. To be more specific, in Sec.~\ref{sec:compsyst} 
we construct an autonomous compactified system on a suitably augmented and extended phase 
space. We then give two simple transformation conditions for the compactified
system to be continuously differentiable on the extended phase space. 
In Sec.~\ref{sec:refenv} we introduce the concept of ``reference envelopes"
to derive  the weakest decay conditions possible
for the existence of the desired compactification transformation.  In Sec.~\ref{sec:sac} 
we distinguish between two-sided compactification for bi-asymptotically constant $\Gamma(t)$ 
and one-sided compactification for asymptotically constant $\Gamma(t)$. 
In Sec.~\ref{sec:compactdyns} we relate the nonautonomous 
dynamics of~\eqref{eq:ode} to the autonomous dynamics of the compactified system. 
Our focus is on the extrapolation of the dynamical structure from the limit systems. This leads us to the construction of stable, unstable and centre invariant manifolds of compact invariant sets 
for the limit systems when embedded in the extended phase space of the compactified system. The geometric shape of these manifolds `encodes' the nonautonomous profile of $\Gamma(t)$, which can greatly simplify analysis of~\eqref{eq:ode}.
In Sec.~\ref{sec:AppA} 
we construct examples of exponential and algebraic compactifications 
that are useful in practice.

\subsection{Motivating Examples}
\label{sec:examples}

Our work is directly motivated by a range of problems 
from applications that can be fitted into the framework of Eq.~\eqref{eq:ode} with 
asymptotically constant $\Gamma(t)$,
and become simpler to analyse after compactification.

\subsubsection{Pullback attractors}
Consider a nonautonomous dynamical system~\eqref{eq:ode} with time $t$ and 
asymptotically constant $\Gamma(t)$ with a past limit
$$
\Gamma(t) \to\Gamma^-\;\;\mbox{as}\;\;t\to -\infty,
$$
and refer to~\cite[Ch.3]{KloedenRasmussen2011} for the notion of a 
 local pullback attractor. 
Each asymptotically stable  compact invariant 
set $A^-$ for the past limit system~\eqref{eq:odea-} can be 
 associated with a local pullback 
attractor of the nonautonomous system~\eqref{eq:ode} as shown in~\cite[Lem.6.2]{Rasmussen2008},~\cite[Th.2.2]{Ashwin2017} 
and~\cite[Th.II.2]{Alkhayuon2018}. When embedded in the extended phase 
space of a compactified system, $A^-$ gains one repelling direction and becomes 
a saddle.  Thus, the process of compactification transforms the local pullback attractor associated with $A^-$ into the unstable or centre invariant 
manifold of $A^-$, which can greatly facilitate analysis of the nonautonomous problem.
Relating pullback attractors to unstable invariant manifolds 
gives an alternative approach to the existing nonautonomous stability theory~\cite{KloedenRasmussen2011}.

\subsubsection{Rate-induced tipping}
A nonlinear nonautonomous dynamical system~\eqref{eq:ode} with  time $t$ 
and bi-asymptotically constant 
$$
\Gamma(t) = \Lambda(rt)\to\lambda^\pm\in\R^d\;\;\mbox{as}\;\;t\to\pm\infty,
$$
is often used to describe nonlinear dynamics of  open systems subject to finite-time external disturbances, growing or decaying external trends, or simply time-varying external inputs $\Lambda(rt)$.  
A rate-induced critical transition,  known as R-tipping~\cite{Ashwin2012}, occurs when the `rate' $r>0$ of the external input exceeds some 
critical value and the system transitions to a different state~\cite{Wieczorek2011,Ashwin2012,Mitry2013,Perryman2014,Ritchie2016,Kiers2018,Alkhayuon2019,OKeeffe2019,Wieczorek2019,Xie2019}.  This genuine nonautonomous bifurcation
is of great interest to natural scientists but cannot, in general, be explained using 
classical bifurcation analysis of the  autonomous frozen system 
$$
\dot{x} = f(x,\lambda),
$$
with a fixed-in-time input parameter $\lambda$.
 The process of compactification transforms the nonautonomous R-tipping
problem into a connecting heteroclinic orbit problem,
which facilitates both numerical
analysis~\cite{Xie2019} and the derivation of rigorous criteria for R-tipping~\cite{Wieczorek2019}. 
Here, the closure of an image of $\Gamma(t)$ has an important meaning: it is
a parameter path in the multi-dimensional parameter space $V\subseteq \R^d$ 
of the frozen system that is traced out by the time-varying external input.

\subsubsection{Radial steady states}
The steady states of nonlinear wave, diffusion or Schr\"odinger equations with a potential ${\mathcal V}(x)$ on $\mathbb{R}^n$ satisfy a semilinear elliptic equation:  
$$
\Delta_x u + f(u) + {\mathcal V}(x)u = 0.
$$
Under the condition that ${\mathcal V}(x)={\mathcal V}(\Vert x\Vert )$, spherically symmetric solutions satisfy the semilinear elliptic boundary value problem (BVP), where 
$r>0$ is the radial direction~\cite{Jones1986,Scheel2003}:
\begin{align}
\label{eq:sphwaves}
u_{rr} + \frac{n-1}{r}\,u_r + f(u) + {\mathcal V}(r)\,u&=0,\;\;
\begin{array}{l}
\lim_{r\to 0}u_r(r)=0,\\
\lim_{r\to +\infty}u(r)=0,
\end{array}
\end{align}
where $u_r=du/dr$. If the potential ${\mathcal V}(r)$ is asymptotically constant with a future 
limit\footnote{This includes the special case ${\mathcal V}(r)= 0$.}, 
the BVP~\eqref{eq:sphwaves} fits naturally into the framework of~\eqref{eq:ode}.
More precisely, by using $r$ as the independent variable and introducing  $v(r)=u_r(r)$ 
as an additional dependent variable so that
$$
x(r)=
\left(
\begin{array}{c}
u(r)\\
v(r)
\end{array}
\right)
\;\;\mbox{and}\;\;\Gamma(r) =  \left(
\begin{array}{c}
\Gamma_1(r)\\
\Gamma_2(r)
\end{array}
\right)=
\left(
\begin{array}{c}
(-n+1)/r\\
-{\mathcal V}(r)
\end{array}
\right)
\to
\left(
\begin{array}{c}
0\\
-{\mathcal V}^+
\end{array}
\right)
\;\;\mbox{as}\;\;r\to +\infty,
$$
 we can rewrite~\eqref{eq:sphwaves} as
\begin{align}
x' = f(x,\Gamma(r)) = \left(
\begin{array}{c}
v\\
\Gamma_1(r)v + \Gamma_2(r)u - f(u)
\end{array}
\right),\;\;
\begin{array}{l}
\lim_{r\to 0} v(r) = 0,\\
\lim_{r\to +\infty}u(r) = 0,
\end{array}
\nonumber
\end{align}
where $'=d/dr$, and $f(0)=0$.

The future limit system corresponds to the 
problem  with constant potential 
and is easily analysed as a Hamiltonian system in the plane~\cite{Jones1986}.
Under common assumptions, the analysis reveals a saddle equilibrium $\eta^+$ at the origin. 
 The full BVP is then solved by finding appropriate trajectories on the stable
or centre-stable invariant manifold of $\eta^+$ in the compactified system.

There are other examples of compactification used for the study of particular solutions of nonlinear equations, such as the lens transform~\cite[Ch.2]{Tao2009}, which is a special case of  pseudoconformal compactification~\cite{Carles2002} for the nonlinear Schr{\"o}dinger equation.

\subsubsection{Stability of nonlinear waves} 

Traveling wave solutions of a reaction-diffusion equation in one space dimension 
\begin{equation}\label{eq:rde}U_t=U_{xx}+F(U),\end{equation}
where $U \in \mathbb{R}^m$, satisfy a second-order system of ordinary differential equations~\cite{Alexander1990}:
 \begin{equation}
 \label{eq:twavesys}
  U''+ c\,U'+F(U)=0,
  \end{equation}
where $'=d/dz$, $z=x-ct$ and $c$ is the speed of the wave. Of interest is commonly a wave $U(z)$ that decays to end-states, $U(z)\to U_\pm$ as $z\rightarrow \pm \infty$. The key in an analysis of the stability of such a wave 
is the eigenvalue problem for the linearisation of (\ref{eq:rde}) at the wave~\cite{Alexander1990}:
\begin{equation}\label{eq:rdelin}
L(U)\,P=P''+ c\,P' + dF(U(z))\,P=\lambda\,P,
\end{equation}
where $dF(U)$ denotes the derivative of $F(U)$ with respect to $U$.
 We can rewrite the second equality in (\ref{eq:rdelin}) as a nonautonomous system on $\mathbb{R}^{2m}$:
 \begin{align}
 \label{eq:twev}
 P'& = Q\\
 Q'&= -c\,Q +\lambda\,P - dF(U(z))\,P,
 \label{eq:twev2}
 \end{align}
with bi-asymptotically constant $dF(U(z))$.
Setting 
 $$
   x(z) = \left( \begin{array}{c} P(z) \\Q(z) \end{array} \right) \;\;\;\mbox{and}\;\;\;
   \Gamma(z) = dF(U(z)) \to dF(U_\pm) \;\;\mbox{as}\;\;z\to \pm\infty,
$$
we can fit~\eqref{eq:twev}--\eqref{eq:twev2} into our general framework
 \begin{align}
x' = f(x,\Gamma(z)) = 
\left(
\begin{array}{c}
Q\\
-c\,Q + \lambda\,P - \Gamma(z)\,P
\end{array}
\right).
\nonumber
\end{align}
 Note that the limit systems here are exactly the eigenvalue problems for the linearisation of~\eqref{eq:rde} at the constant (time and space independent) solutions $U_\pm$. 
 
 The Evans Function is the primary tool used to capture these eigenvalues. The construction 
of it in~\cite{Alexander1990}, where the first general formulation was given, 
involves a two-sided compactification of system~\eqref{eq:twev}--\eqref{eq:twev2}.

\section{ Compactification}
\label{sec:compact}

\begin{figure}[t]
  \begin{center}
    \includegraphics[width=15.cm]{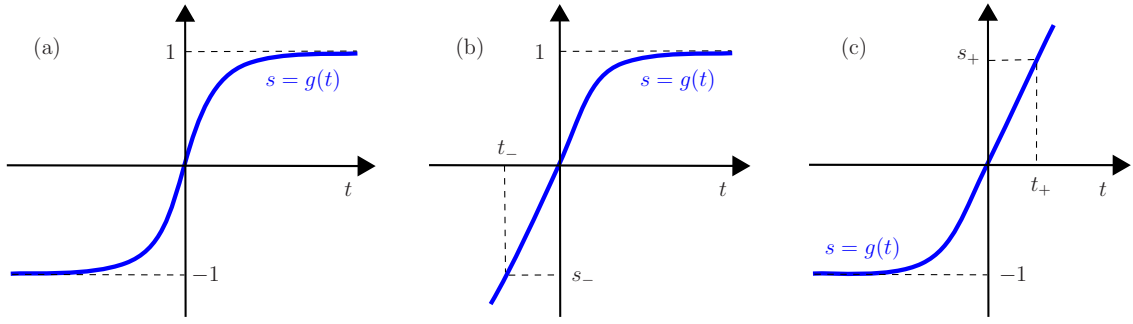}
  \end{center}
  \vspace{-3mm}
  \caption{ Examples of the coordinate transformation $s=g(t)$
  from (a) Assumption~\ref{asmp:g}, (b) Assumption~\ref{asmp:g+} and
  (c) Assumption~\ref{asmp:g-}.
 }
  \label{fig:g}
\end{figure}

The  aim of this section is to reformulate the nonautonomous 
system~\eqref{eq:ode} into an autonomous system so that:
\begin{itemize}
    \item 
    The new system contains the autonomous flow and compact invariant sets of the future limit system~\eqref{eq:odea+} and/or the past limit system~\eqref{eq:odea-}.
    \item
    The dimension of the vector field increases just by one, independently of the dimension 
    and monotonicity of the nonautonomous term $\Gamma(t)$.
\end{itemize}
We show that this can be achieved 
under quite general assumptions on ${\Gamma}(t)$ by augmenting system~\eqref{eq:ode} with 
$s = g(t)$ depicted in Fig.~\ref{fig:g} as an additional dependent variable. In other words, 
instead of having a problem with the additional dimension being unbounded as in the usual 
extended system~\eqref{eq:odeext0}, we augment system~\eqref{eq:ode} so that the additional 
dimension becomes a compact interval. Specifically, we:
\begin{itemize}
 \item[(i)] 
 Restrict to asymptotically and bi-asymptotically constant $\Gamma(t)$.
  \item[(ii)] 
 Compactify the real $\nu$-line into the compact 
  $s$-interval $[-1,1]$  if $\Gamma(t)$ is bi-asymptotically constant.
  \item[(iii)]
Compactify the half $\nu$-line  $[\nu_-,+\infty)$  into the compact   $s$-interval $[s_-,1]$  if $\Gamma(t)$ is 
asymptotically constant with a future limit, or compactify the half $\nu$-line $(-\infty,\nu_+]$ into the compact  
$s$-interval $[-1,s_+]$ if $\Gamma(t)$ is asymptotically constant with a past limit.
\end{itemize}
%

 \subsection{Two-Sided Compactification for Bi-Asymptotically Constant $\Gamma(t)$}
\label{sec:twopointcomp}

We reformulate the nonautonomous system~\eqref{eq:ode}  with a bi-asymptotically constant 
$\Gamma(t)$ into a compactified system that is autonomous and contains the flow
and compact invariant sets of the 
future~\eqref{eq:odea+} and past~\eqref{eq:odea-} limit systems. This is achieved via two-sided 
compactification that uses  the coordinate transformation $s=g(t)$ depicted in Fig.~\ref{fig:g}(a).
More precisely, we assume
\begin{asmp}
\label{asmp:g}
A coordinate transformation $s=g(t)$ maps the real $t$-line  onto the finite  $s$-interval $(-1,1)$, is at lest $C^2$-smooth, bi-asymptotically constant with future limit 1 and past limit $-1$, and strictly increasing with vanishing first derivative  as $t$ tends to $\pm\infty$:
\begin{equation}
\label{eq:g}
g:\R\to(-1,1),\;\; g\in C^{k\ge 2},\;\;
\lim_{t\to \pm\infty}g(t) = \pm 1,\;\;
\dot{g}(t)>0\;\;\mbox{for}\;\;t \in\R\;\;\mbox{and}\;\;\lim_{t\to \pm\infty}\dot{g}(t) = 0.
\end{equation} 
\end{asmp}
\begin{rmk}
In practice, we introduce  {\em compactification parameter(s)} to control the rate/order of asymptotic decay of $g(t)$, and work with {\em parametrised 
compactifications}; see Sec.~\ref{sec:AppA} ahead. 
The dependence on the compactification parameter(s) does not affect, and is thus 
left out of, general statements of this section.
\end{rmk}

\subsubsection{Compactified System: Transformation Conditions}
\label{sec:compsyst}

Compactification is a three-step process. The first step is to
make the additional dimension bounded by augmenting the asymptotically 
autonomous system~(\ref{eq:ode}) with $s=g(t)$ as an additional dependent variable
\begin{align}
  \label{eq:odeexta}
  \dot{x}&= f\left(x,\Gamma\left(g^{-1}(s)\right)\right),\\
  \label{eq:odeext2a}
  \dot{s}&=  \dot{g}\left(g^{-1}(s)\right),\quad\quad\quad\; s(t_0)=g(t_0).
\end{align}
Since the practical implementation of the compactification
requires the inverse coordinate transformation $t=h(s):= g^{-1}(s)$,
we reformulate~\eqref{eq:odeexta}--\eqref{eq:odeext2a} in terms of $h(s)$ alone:
\begin{align}
  \label{eq:odeextb}
  \dot{x}&= f(x,\Gamma(h(s))),\\
  \label{eq:odeext2b}
  \dot{s}&=  \gamma(s),
\end{align}
where $\gamma(s)=\dot{g}\left(h(s)\right) = 1/h'(s)$ is the augmented component of the vector
field, $\dot{} = d/dt$, and $'=d/ds$.
System~\eqref{eq:odeextb}--\eqref{eq:odeext2b} is defined on
$U\times(-1,-1)$, which is an open subset of $\R^{n+1}$.
The second step is to make the $s$-interval closed
by including $s=\pm 1$ ($t=\pm\infty$) and extending the augmented by~\eqref{eq:odeext2b} 
vector field to subspaces $\{s=\pm 1\}$:
\begin{align}
  \label{eq:lambda_s}
  f(x,\Gamma(h(s)))=&
  \left\{
    \begin{array}{rcl}
      f(x,\Gamma(h(s))) &\mbox{for}& s\in(-1,1),\\
      f(x,\Gamma^-) &\mbox{for}& s=-1,\\
      f(x,\Gamma^+) &\mbox{for}& s=1,
    \end{array}
  \right.\\
  \label{eq:lambda_q}
  \gamma(s)=&
  \left\{
    \begin{array}{rcl}
       1/h'(s) &\mbox{for}& s\in(-1,1),\\
      0 &\mbox{for}& s=\pm 1.
    \end{array}
  \right.
\end{align}
This gives an autonomous {\em compactified system}~\eqref{eq:odeextb}--\eqref{eq:lambda_q}
that is defined on the extended phase space $U\times[-1,1]$. Most importantly, subspaces 
$\{s=1\}$ and $\{s=-1\}$ are flow-invariant and carry the autonomous dynamics and compact invariant sets 
of the future~\eqref{eq:odea+} and past~\eqref{eq:odea-} limit systems, respectively. However, 
it is not generally guaranteed that the extended vector field is differentiable at the added 
invariant subspaces $\{s=\pm 1\}$. The third step is to give testable criteria for 
the extended vector field to be continuously differentiable.
Examining first-order partial derivatives of the extended vector field~\eqref{eq:lambda_s}--\eqref{eq:lambda_q}
leads to the following result:
\begin{proposition}
\label{l2}
  {\em (Transformation Conditions.)}
 Consider a nonautonomous system~(\ref{eq:ode}) on $U$ with $C^1$-smooth $f$ and $\Gamma$, and 
 bi-asymptotically constant $\Gamma(t)$. For a chosen  coordinate transformation $g(t)$ from 
 Assumption~\ref{asmp:g}, with the inverse $h(s) = g^{-1}(s)$, the ensuing compactified system~\eqref{eq:odeextb}--\eqref{eq:lambda_q}  
 is $C^{1}$-smooth on $U\times[-1,1]$ if and only if
 \begin{align}
  \label{eq:dLdg}
      \lim_{t\to \pm\infty}\,&\frac{\dot{\Gamma}(t)}{\dot{g}(t)} = \lim_{s\to \pm1^\mp}\frac{d}{ds}\Gamma(h(s))\;\;\mbox{exist},\\
    \label{eq:d2gdt2}
      \lim_{t\to \pm\infty}\,&\frac{\ddot{g}(t)}{\dot{g}(t)}  = -\lim_{s\to
        \pm1^\mp}\frac{h''(s)}{\left(h'(s)\right)^2}\;\;\mbox{exist}.
    \end{align}
\end{proposition}
%
\begin{figure}[t]
  \begin{center}
    \includegraphics[width=13.cm]{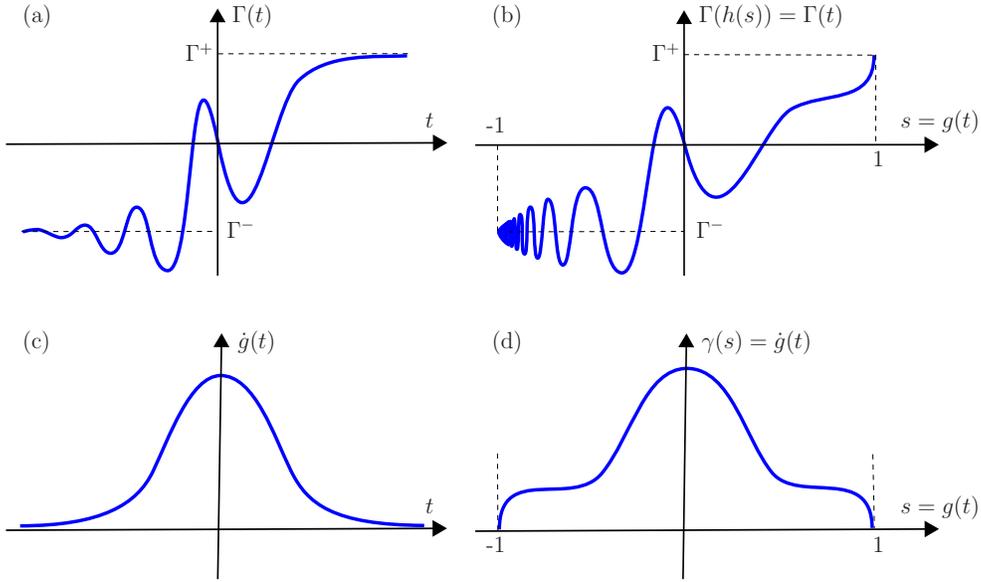}
  \end{center}
  \vspace{-3mm}
  \caption{(a)--(b) Example of transformation-induced loss of differentiability of 
    $\Gamma(h(s))$ at $s=\pm 1$ when condition~(\ref{eq:dLdg}) is violated.
    (c)--(d) Example of $\gamma(s)$ that is 
    non-differentiable at $s=\pm 1$ when
    condition~(\ref{eq:d2gdt2}) is violated.  
  }
  \label{fig:comp_cond}
\end{figure}
%
The one-dimensional examples from Fig.~\ref{fig:comp_cond} give insight 
into transformation conditions~(\ref{eq:dLdg}) and~(\ref{eq:d2gdt2}),
which ensure differentiability of $\Gamma(h(s))$ and $\gamma(s)$, respectively, 
at $s=\pm 1$. 
Recall that $\dot{g}(t)$ limits to zero as $t\to\pm\infty$ and note the following.
The slope of $\Gamma(h(s))$ in Fig.~\ref{fig:comp_cond}(b) 
is given by 
$d\Gamma/ds = d\Gamma/dg = \dot{\Gamma}/\dot{g}$. Thus, the derivative of
$\Gamma(h(s))$ becomes undefined  at $s=\pm1$ if $\dot{\Gamma}(t)$
does not limit to zero or if $\dot{\Gamma}(t)$ approaches zero 
slower than $\dot{g}(t)$ as $t\to\pm\infty$; see Sec.~\ref{sec:refenv} for optimal conditions on the asymptotic decay of $\dot{\Gamma}(t)$.
In other words, coordinate transformations
that decay faster than the 
nonautonomous term 
violate condition~(\ref{eq:dLdg}) and cause 
transformation-induced 
loss of differentiability of $\Gamma(h(s))$ at $s=\pm 1$. 
Similarly, the slope of $\gamma(s)$ in Fig.~\ref{fig:comp_cond}(d)
is given by $d\gamma/ds = d\dot{g}/ds = d\dot{g}/dg = \ddot{g}/\dot{g}$. Thus, 
the derivative of $\gamma(s)$ becomes undefined at $s=\pm 1$ if 
$\ddot{g}(t)$ does not limit to zero or if $\ddot{g}(t)$
approaches zero slower than $\dot{g}(t)$ as $t\to \pm\infty$.
In other words, 
coordinate transformations whose speed tends to zero faster than acceleration violate condition~(\ref{eq:d2gdt2}) and make 
$\gamma(s)$ non-differentiable at $s=\pm 1$.
The following remark gives additional insight into 
condition~(\ref{eq:d2gdt2}).
\begin{rmk}
\label{rmk:superexp}
Transformation condition~(\ref{eq:d2gdt2}) excludes coordinate transformations 
$g(t)$ that decay faster than exponentially.
To see that, consider
\begin{displaymath}
g(t) \sim \left\{
 \begin{array}{r c l}
    1 - \exp\left(-t^k\right)& \mbox{as} & t\to +\infty,\\
-1 + \exp\left(-(-t)^k\right)& \mbox{as} & t\to -\infty,
  \end{array}
  \right. 
\end{displaymath}
for $k>0$, and note that
\begin{displaymath}
\lim_{t\to \pm\infty}\frac{\ddot{g}(t)}{\dot{g}(t)}= 
\mp k
\lim_{t\to \pm\infty} (\pm t)^{k-1}= \left\{
  \begin{array}{r c l}
    \mp\infty & \mbox{if} & k > 1,\\
    \mp 1 & \mbox{if} & k=1,\\
    0 & \mbox{if} &  0 < k < 1.
  \end{array}
\right.
\end{displaymath}
 This can  be understood intuitively via behaviour near invariant subspaces: 
solutions of ODEs cannot approach invariant subspaces faster than exponentially unless 
they blow up and cease to exist.
\end{rmk}

\noindent
{\em Proof of Proposition~\ref{l2}}\\
$C^1$-smoothness of the extended vector field~\eqref{eq:lambda_s}--\eqref{eq:lambda_q}
on $U\times(-1,1)$ follows from $C^1$-smoothness of $f$ on $U$, 
$C^1$-smoothness of $\Gamma$ on $\R$,  and $C^2$-smoothness of $g$ on $\R$.

What needs to be examined is $C^1$-smoothness of the extended vector field at the added
right and left invariant subspaces $\{s = \pm 1\}$.
It follows from Definition~\ref{defn:ac} and Assumption~\ref{asmp:g} 
that $\Gamma(h(s))$ and $\gamma(s)$ are continuous at $s=\pm 1$:
\begin{align}
\label{eq:st}
&\lim_{s\to\pm 1^\mp} h(s) = \lim_{t\to\pm \infty} t,\\
&\lim_{s\to\pm 1^\mp} \Gamma(h(s)) = \lim_{t\to\pm \infty} \Gamma(t) = \Gamma^\pm = \Gamma(h(\pm 1)),\nonumber\\
&\lim_{s\to\pm 1^\mp} \gamma(s) = \lim_{s\to\pm 1^\mp} \dot{g}\left(h(s)\right) = \lim_{t\to\pm \infty}\dot{g}(t)  = 0 = \gamma(\pm 1).\nonumber
\end{align}
Thus, the extended vector field is continuous on $U\times[-1,1]$.
The first derivative of the extended vector field  is 
continuous at $\{s=\pm 1\}$ if the left- and 
right-sided limits $s\to \pm1^\mp$ exist for all first-order partial derivatives.
To check this, consider the Jacobian of
the extended vector field
\begin{equation}
\label{eq:jac}
J(x,s)=
\begin{pmatrix}
  \left(\frac{\partial f}{\partial x}\right)_{n\times n} &
    \left(\frac{\partial f}{\partial s}\right)_{n\times
    1} 
\vspace*{3mm}
\\
  (0)_{1\times n} & \left(\frac{d\gamma}{ds}\right)_{1\times 1}
\end{pmatrix},
\end{equation}
where the subscripts indicate the size of the matrix components of $J(x,s)$, and 
\begin{align}
\label{eq:dfds}
\left(\frac{\partial f}{\partial s}\right)_{n\times
    1} &=
\left(\frac{\partial f}{\partial\Gamma}\right)_{n\times d}
\left(\frac{d}{d s}\Gamma(h(s))\right)_{d\times 1}.
\end{align}
Using the chain rule and the inverse function theorem gives
\begin{align}
\label{eq:l2a}
\frac{d}{d s}\Gamma(h(s)) = 
\dot{\Gamma}(h(s))\,h'(s) =
\frac{\dot{\Gamma}(h(s))}{\dot{g}\left(h(s)\right)}.
\end{align}
Noting that $h(s)$ is at least $C^2$-smooth by the (general) inverse function theorem, 
 differentiate $s = g\left(h(s) \right)$ twice with respect to $s$ to obtain
\begin{equation}
  \label{eq:d2hds2n}
\ddot{g}(h(s))\,\left(h'(s)\right)^2 = - \dot{g}(h(s))\, h''(s),
\end{equation}
which, together with the chain rule and the fact that $h'(s) = 1/\dot{g}(h(s))>0$, gives
\begin{align}
\label{eq:l2b}
\gamma'(s)&:=\frac{d}{ds}\gamma(s) = \ddot{g}\left(h(s)\right)\,h'(s)
= \frac{\ddot{g}(h(s))}{\dot{g}(h(s))} 
= - \frac{h''(s)}{\left(h'(s)\right)^2}.
\end{align}
It follows from the continuity of $\Gamma(h(s))$ at $s=\pm 1$ that 
 $(\partial f/\partial x)_{n\times n}$ and $(\partial
f/\partial \Gamma)_{n\times d}$ are continuous at $s=\pm 1$.
It then follows from~\eqref{eq:st},~\eqref{eq:dfds} and~\eqref{eq:l2a} that 
$\left(\partial f/\partial s \right)_{n\times 1}$
is continuous at $s=\pm 1$ if and only if the first transformation condition~\eqref{eq:dLdg} 
is satisfied.  
It follows from~\eqref{eq:st} and~\eqref{eq:l2b} that the $J_{n+1,n+1}$ component of the Jacobian is continuous 
at $s=\pm 1$ if and only if the second transformation condition~\eqref{eq:d2gdt2} is satisfied. 
Thus, the first derivative of the extended vector field is continuous on $U\times[-1,1]$ 
if and only if transformation conditions~\eqref{eq:dLdg} and~\eqref{eq:d2gdt2} are satisfied.
\qed\\

\subsubsection{Compactified System: Optimal Existence Criteria}
\label{sec:refenv}

Here we derive optimal criteria on $\dot{\Gamma}(t)$ that guarantee a 
suitable coordinate transformation $g(t)$ that satisfies conditions~\eqref{eq:dLdg} and~\eqref{eq:d2gdt2} from Proposition~\ref{l2} exists. 
The derivation of the existence criteria is guided by the observation that `normal' examples of bi-asymptotically constant $\Gamma(t)$
have two additional and desirable properties. Firstly,  $\dot{\Gamma}(t)$ limits to 
zero as $t$ tends to positive and negative infinity. Secondly, the asymptotic approach 
of $\dot{\Gamma}(t)$ towards zero is not slower than the asymptotic approach of $\Gamma(t)$ 
towards $\Gamma^\pm$.

\begin{figure}[t]
  \begin{center}
    \includegraphics[width=8.cm]{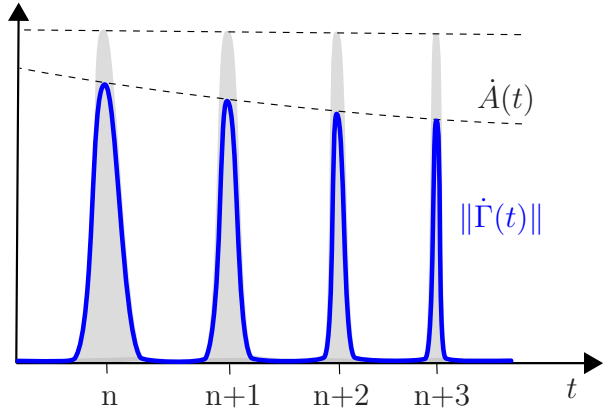}
  \end{center}
  \vspace{-3mm}
  \caption{
    Consider  $C^{1}$-smooth bumps of equal amplitudes whose (shaded) areas form a convergent
    series as $t\to +\infty$.  Then, $\Vert \dot{\Gamma}(t)\Vert$ is set to consist of
    such bumps with bump amplitudes remaining constant or decaying arbitrarily slowly to zero as $t\to +\infty$. Since the total area below $\Vert\dot{\Gamma}(t)\Vert$ is finite by construction, we have $\Gamma(t)\to \Gamma^+\in\R$ as $t\to +\infty$ in 
    spite of the envelope $\dot{A}(t)$ remaining constant or decaying arbitrarily 
    slowly to zero.
  }
  \label{fig:Latypical}
\end{figure}

However, $\Gamma(t)\to \Gamma^\pm\in\R^d$ does not imply $\dot{\Gamma}(t)\to 0$ in
general, meaning that there exist `pathological' examples of $\Gamma(t)\to
\Gamma^\pm$ whose  derivatives do not have a future or past limit, or
approach zero 
arbitrarily slowly as $t$ tends to positive or negative infinity.  One example are 
damped oscillations with increasing frequency, where the frequency increase `beats' 
the amplitude decay. For example, $\Gamma(t)\sim \sin(t^2)/t$ as $t\to+\infty$ has a future 
limit but its first derivative $\dot{\Gamma}(t) \sim 2\cos(t^2) -\sin(t^2)/t^2$ does not. 
Another example is depicted in Fig.~\ref{fig:Latypical} where $\Gamma(t)\to \Gamma^+$ 
but $\dot{\Gamma}(t)$ may not have a future limit or may approach zero arbitrarily 
slowly.
Conversely, $\dot{\Gamma}(t)\to 0$ as $t\to\pm\infty$  does not imply 
$\Gamma(t)\to \Gamma^\pm\in\R^d$ either.  For example,  $\Gamma(t)\sim\ln(t)$ as 
$t\to +\infty$ does not have a future limit even though its first derivative 
$\dot{\Gamma}(t)\sim 1/t$ limits to zero.

To exclude `pathological' examples of bi-asymptotically constant $\Gamma(t)$ it is convenient to work with a monotone envelope of $\Vert\dot{\Gamma}(t)\Vert$, denoted
by $\dot{A}(t)$ as depicted in Fig.~\ref{fig:Latypical}, and ask about the slowest-decaying $\dot{A}(t)\to 0$ so that $A(t)$ has a future and past limits.  Formulating this question in terms 
of the integral
\begin{align}
\label{eq:env}
\lim_{t\to +\infty} A(t) = \int_{t_0}^\infty
\dot{A}(\mu)\,d\mu= A^+\in\R,
\end{align}
shows that such a slowest-decaying envelope does not
exist: given any $\dot{A}(t)\to 0$ that
satisfies~(\ref{eq:env}) for some $t_0\in\R$, one can construct a slower-decaying one that
also satisfies~(\ref{eq:env}).  Nonetheless, it is possible to work with
a parametrised family of envelopes that satisfy~(\ref{eq:env}) and
can be chosen to decay sufficiently slowly for the problem at hand. Specifically, we consider
\begin{align}
\label{eq:A}
A(t,m) = -\frac{1}{\ln^m(|t|)},\nonumber
\end{align}
together with its derivative
\begin{equation}
\label{eq:Lmeps+}
\dot{A}(t,m) =
\left(t\,\ln^m(|t|)\prod_{k=1}^m\ln^k(|t|)\right)^{-1},
\end{equation}
parametrised by a non-negative integer $m\in\N_0$, where
$|t| > \exp^{m-2}(e)$ and $\ln^m$ denotes a composition 
of $m$ logarithmic functions.
\begin{defn}
\label{def:refenv}
We call $\dot{A}(t,m)$ from Eq.~\eqref{eq:Lmeps+} the {\em parametrised reference envelope}.
\end{defn}
The reference envelope formula~\eqref{eq:Lmeps+} looks rather technical but one can  gain further  
insight by  writing out examples of $A$ and $\dot{A}$ with  $m=2$ that  are defined for  $|t| > e$:
\begin{equation}
\label{eq:Lmeps+2}
A(t,2) =-\frac{1}{\ln(\ln(|t|)}\;\;\;{\rm and}\;\;\;
\dot{A}(t,2) =\frac{1}{t\,\ln(|t|)\,\left(\ln(\ln(|t|))\right)^{2}}.\nonumber
\end{equation}
More generally, one can verify the following properties
\begin{equation}
\label{eq:Aprop}
\lim_{t\to \pm\infty}{A}(t,m)= 0,\;\;
\dot{A}(t,m)
\left\{
\begin{array}{r c l}
> 0 &\mbox{for} & t > \exp^{m-2}(e),\\
< 0 &\mbox{for} & t < -\exp^{m-2}(e),
\end{array}
\right.
\,\;\mbox{and}\;\;\lim_{t\to \pm\infty}\dot{A}(t,m)=0.
\end{equation}

We now use the reference envelope concept to restrict to `normal'  examples of bi-asymptotically constant  $\Gamma(t)$ and ensure the existence of a continuously differentiable compactified system.
\begin{defn}
\label{def:strongac}
We call a bi-asymptotically constant $\Gamma(t)$ {\em normal} if there is an $m$ such that
\begin{align}
\label{eq:ac2}
&\lim_{t\to \pm\infty} \frac{\dot{\Gamma}(t)}{\dot{A}( t,m)}\;\;\mbox{exist}.
\end{align}
\end{defn}
\begin{theorem}
\label{th2}
{\em (Existence of a $C^1$-smooth Compactified System.)}
Consider a nonautonomous system~(\ref{eq:ode}) with with $C^1$-smooth $f$ and $\Gamma$.
If $\Gamma(t)$ is normal bi-asymptotically constant then there exists 
a coordinate transformation $s=g(t)$ from Assumption~\ref{asmp:g} such that
the extended vector field~(\ref{eq:lambda_s})--(\ref{eq:lambda_q}) of the compactified 
system is  $C^{1}$-smooth on $U\times[-1,1]$.
  \end{theorem}
In other words, the normal bi-asymptotic constant condition~\eqref{eq:ac2} gives  sufficient existence criteria that are both 
{\em optimal} and {\em easily testable}. 
To be more precise, the criteria are optimal in the sense that they 
identify the weakest decay of $\dot{\Gamma}(t)$ possible to eliminate `pathological' bi-asymptotically constant $\Gamma(t)$ and enable construction of a
continuously differentiable
compactified system, while still allowing super-exponential, 
exponential, algebraic, logarithmic or even sub-logarithmic decay of $\Gamma(t)$. 
\vspace{3mm}

\noindent
{\em Proof of Theorem~\ref{th2}}\\
\noindent
Use the parametrised reference envelope from Def.~\ref{def:refenv} to construct 
a parametrised by $m$ coordinate transformation $g_{(m)}:\R\to(-1,1)$ that is at least 
$C^2$-smooth, strictly increasing and
\begin{align}
g_{(m)}(t) \sim \left\{
 \begin{array}{l l c}
     1 + A(t,m)  &
     \mbox{as} & t\to+\infty,\\
     -1 - A(t,m) &
     \mbox{as} & t\to -\infty.
     \end{array}\nonumber
\right.
\end{align}
It follows from the asymptotic properties of  $A(t,m)$  that
\begin{align}
\lim_{t\to\pm\infty} g_{(m)}(t)  = \pm 1\;\;\mbox{and}\;\;
\lim_{t\to\pm\infty} \dot{g}_{(m)}(t)  = 0.\nonumber
\end{align}
Thus, $g_{(m)}(t)$ satisfies Assumption~\ref{asmp:g} by construction.
By Definition~\ref{def:strongac}, for any normal bi-asymptotically constant $\Gamma(t)$ 
there is an $m$ such that
$$
\lim_{t\to \pm\infty}\frac{\dot{\Gamma}(t)}{\dot{g}_{(m)}(t)} = 
\pm  \lim_{t\to \pm\infty}\frac{\dot{\Gamma}(t)}{\dot{A}( t,m)},
$$
exist, meaning that $g_{(m)}(t)$ satisfies the first transformation 
condition~\eqref{eq:dLdg} by construction. Finally, one can verify by 
induction that 
$$
\lim_{t\to \pm\infty}\frac{\ddot{g}_{(m)}(t)}{\dot{g}_{(m)}(t)} = \lim_{t\to \pm\infty}\frac{\ddot{A}(t,m)}{\dot{A}(t,m)}=0,
$$
for any non-negative $m\in\N_0$, meaning that $g_{(m)}(t)$ satisfies 
the second transformation condition~\eqref{eq:d2gdt2} by construction.
It then follows from Proposition~\ref{l2} that the extended vector field~(\ref{eq:lambda_s})--(\ref{eq:lambda_q}) is  $C^{1}$-smooth on $U\times[-1,1]$.
 \qed

\subsection{One-Sided Compactification for Asymptotically Constant $\Gamma(t)$}
\label{sec:sac}

Here we discuss briefly two one-sided subcases of the two-sided compactification from Sec~\ref{sec:twopointcomp}.
The discussion does not require new analysis, but may be helpful to readers interested in problems 
with asymptotically constant $\Gamma(t)$.
For asymptotically constant $\Gamma(t)$ with a future limit, we reformulate the nonautonomous system~\eqref{eq:ode} into a compactified system that is autonomous and contains the flow 
and compact invariant sets of the future limit system~\eqref{eq:odea+}. This is achieved 
via one-sided compactification that uses coordinate transformation $s=g(t)$ depicted in Fig.~\ref{fig:g}(b). For asymptotically constant $\Gamma(t)$ with a past limit, we reformulate 
the nonautonomous system~\eqref{eq:ode} into a compactified system that is autonomous and 
contains the flow and compact invariant sets of the past limit system~\eqref{eq:odea-}. 
This is achieved via one-sided compactification that uses coordinate transformation $s=g(t)$ 
depicted in Fig.~\ref{fig:g}(c).

Similar to Def.~\ref{def:strongac}, we restrict to `normal'  examples of asymptotically constant  $\Gamma(t)$:
\begin{defn}
\label{def:strongac2}
We call an asymptotically constant $\Gamma(t)$ with a future limit {\em normal} if there is an $m$ such that the $t\to +\infty$ limit in~\eqref{eq:ac2} exists.
We call an asymptotically constant $\Gamma(t)$ with a past limit {\em normal} if there is an $m$ such that the $t\to -\infty$ limit in~\eqref{eq:ac2} exists.
\end{defn}

\subsubsection{Right-Sided Compactification}
Consider nonautonomous system~\eqref{eq:ode} with asymptotically constant $\Gamma(t)$ with  
future limit $\Gamma^+$, and assume that
\begin{asmp}
\label{asmp:g+}
A coordinate transformation $s=g(t)$ maps the half $t$-line $[t_-,+\infty)$ onto the finite  $s$-interval $[s_-,1)$, is at least $C^2$-smooth, asymptotically constant with future limit $1$, 
and strictly increasing with vanishing first derivative as $t\to +\infty$:
\begin{equation}
\label{eq:g+}
g:[t_-,+\infty)\to[s_-,1),\;\; g\in C^{k\ge 2},\;\;
\lim_{t\to +\infty}g(t) = 1,\;\;\dot{g}(t)>0\;\;\mbox{for}\;\;t \ge t_-\;\;\mbox{and}\;\;\lim_{t\to +\infty}\dot{g}(t) = 0.
\end{equation} 
\end{asmp}
\noindent
Then consider the following autonomous {\em compactified system}
\begin{align}
  \label{eq:odeext}
  \dot{x}&= f(x,\Gamma(h(s))),\\
  \label{eq:odeext2}
  \dot{s}&=  \gamma(s),\\
  \label{eq:lambda_s+}
  f(x,\Gamma(h(s)))=&
  \left\{
    \begin{array}{rcl}
      f(x,\Gamma(h(s))) &\mbox{for}& s\in[s_-,1),\\
      f(x,\Gamma^+) &\mbox{for}& s=1,
    \end{array}
  \right.\\
  \label{eq:lambda_q+}
  \gamma(s)=&
  \left\{
    \begin{array}{rcl}
      1/h'(s) &\mbox{for}& s\in[s_-,1),\\
      0 &\mbox{for}& s=1,
    \end{array}
  \right.
\end{align}
that is defined on the extended phase space $U\times[s_-,1]$ with
flow-invariant subspace $\{s=1\}$ that carries the autonomous dynamics and compact 
invariant sets of the future limit system~\eqref{eq:odea+}.

Proposition~\ref{l2} and Theorem~\ref{th2}
apply to compactified system~\eqref{eq:odeext}--\eqref{eq:lambda_q+} after 
we leave out the limits $t\to -\infty$ ($s\to -1^+$) and replace:
``bi-asymptotically constant" with ``asymptotically constant with a future limit",
``compactified system~\eqref{eq:odeextb}--\eqref{eq:lambda_q}" with ``compactified system~\eqref{eq:odeext}--\eqref{eq:lambda_q+}",
``phase space $U\times[-1,1]$" with ``phase space $U\times[s_-,1]$" and
``Assumption~\ref{asmp:g}" with ``Assumption~\ref{asmp:g+}".

\subsubsection{Left-Sided Compactification}
Similarly, consider nonautonomous system~\eqref{eq:ode} with asymptotically constant $\Gamma(t)$ with  
past limit $\Gamma^-$, and assume that
\begin{asmp}
\label{asmp:g-}
A coordinate transformation $s=g(t)$ maps the half $t$-line $(-\infty,t_+]$ onto the finite  $s$-interval $(-1,s_+]$, is at least $C^2$-smooth, asymptotically constant with past limit $-1$, and strictly increasing with vanishing first derivative as $t\to -\infty$:
\begin{equation}
\label{eq:g-}
g:(-\infty,t_+]\to(-1,s_+],\;\; g\in C^{k\ge 2},\;\;
\lim_{t\to -\infty}g(t) = -1,\;\;\dot{g}(t)>0\;\;\mbox{for}\;\;t \le t_+
\;\;\mbox{and}\;\;\lim_{t\to -\infty}\dot{g}(t) = 0.
\end{equation} 
\end{asmp}
\noindent
Then consider the following autonomous {\em compactified system}
\begin{align}
  \label{eq:odeext-}
  \dot{x}&= f(x,\Gamma(h(s))),\\
  \label{eq:odeext2-}
  \dot{s}&=  \gamma(s),\\
  \label{eq:lambda_s-}
 f(x,\Gamma(h(s)))=&
  \left\{
    \begin{array}{rcl}
      f(x,\Gamma(h(s))) &\mbox{for}& s\in(-1,s_+],\\
      f(x,\Gamma^-) &\mbox{for}& s=-1,
    \end{array}
  \right.\\
  \label{eq:lambda_q-}
  \gamma(s)=&
  \left\{
    \begin{array}{rcl}
      1/h'(s) &\mbox{for}& s\in(-1,s_+],\\
      0 &\mbox{for}& s=-1,
    \end{array}
  \right.
\end{align}
that is defined on the extended phase space $U\times[-1,s_+]$ with
flow-invariant subspace $\{s=-1\}$ that carries the autonomous dynamics and compact 
invariant sets of the past limit system~\eqref{eq:odea-}.

Proposition~\ref{l2} and Theorem~\ref{th2}
apply to compactified system~\eqref{eq:odeext-}--\eqref{eq:lambda_q-} after 
we leave out the limits $t\to +\infty$ ($s\to 1^-$) and replace:
``bi-asymptotically constant" with ``asymptotically constant with a past limit",
``compactified system~\eqref{eq:odeextb}--\eqref{eq:lambda_q}" with ``compactified 
system~\eqref{eq:odeext-}--\eqref{eq:lambda_q-}",
``phase space $U\times[-1,1]$" with ``phase space $U\times[-1,s_+]$" and
``Assumption~\ref{asmp:g}" with ``Assumption~\ref{asmp:g-}".

\section{Compactified System Dynamics}
\label{sec:compactdyns}

In this section we relate the nonautonomous dynamics of system~\eqref{eq:ode} with 
bi-asymptotically constant $\Gamma(t)$ and the  autonomous dynamics of the compactified 
system~\eqref{eq:odeextb}--\eqref{eq:lambda_q}.
Specifically, we focus on the extrapolation of dynamical structure 
from (one of) the limit systems~\eqref{eq:odea+} or~\eqref{eq:odea-}, 
start with a general remark, and follow with rigorous statements
that can greatly simplify analysis of~\eqref{eq:ode}.
Here, we distinguish between equilibria and more general 
compact invariant sets:
\vspace*{-2mm}
\begin{itemize}
    \item $A$ denotes a general compact invariant set.
    \item $\eta$ denotes a compact invariant set that is an equilibrium point.
\end{itemize}

\begin{rmk}
\label{rmk:csd}
In the compactified system~\eqref{eq:odeextb}--\eqref{eq:lambda_q},
the time evolution of $s(t)$ is not influenced by the time evolution of $x(t)$.
Owing to this special skew-product structure:
\begin{itemize}
\item[(i)]
Any  compact invariant set {$A^+$} for the future limit 
system~\eqref{eq:odea+} gains one attracting direction, and any 
compact invariant set {$A^-$} for the past limit system~\eqref{eq:odea-} 
gains one repelling direction when embedded in the extended phase space 
of~\eqref{eq:odeextb}--\eqref{eq:lambda_q}, 
in the sense  of monotone increasing $s(t)$:
$$
\dot{s}(t)=\dot{g}(t) >0.
$$
\item[(ii)]
A regular~\footnote{A compact invariant set is {\em regular} if the full Lyapunov spectrum exists.} compact invariant set {$A^\pm$} of one of the limit systems can be viewed as gaining an additional Lyapunov exponent~\footnote{Referred to as the {\em normal Lyapunov exponent} in the terminology of~\cite{Ashwin1996}.}, denoted $l_s^\pm$,  when embedded in the extended phase space. 
The additional Lyapunov exponent 
quantifies linear stability  in the $s$-direction, is independent of $x$, and is given by the 
$(n+1,n+1)$-th element of the Jacobian~\eqref{eq:jac} at $\{s=\pm 1\}$,
or by the limit in the second transformation condition~\eqref{eq:d2gdt2}:
\begin{equation}
\label{eq:sLyap}
l^\pm_{s} = 
\gamma'(s)\bigg|_{s=\pm 1} = 
- \frac{h''(s)}{\left(h'(s)\right)^2}\bigg|_{s=\pm 1} =
\lim_{t\to \pm\infty}\,\frac{\ddot{g}_{}(t)}{\dot{g}_{}(t)},
\end{equation}
where the second equality follows from Proposition~\ref{l2}.
We note that $l^\pm_{s}$ is zero when $g(t)$ decays slower than exponentially; 
see Sec.~\ref{sec:ac}.
\item[(iii)]
The Lyapunov vector corresponding to the additional Lyapunov exponent is independent of $x$. The vector is normal to $\{s=\pm 1\}$ if the top $n$ elements
in the last column of the Jacobian~\eqref{eq:jac} at $\{s=\pm 1\}$ are zero, or if the first transformation condition~\eqref{eq:dLdg} is zero
$$
\frac{\partial f}{\partial s} \bigg|_{s=\pm 1}=
\frac{\partial f}{\partial \Gamma}\bigg|_{s=\pm 1}\, 
\frac{d}{ds}\Gamma(h(s)) \bigg|_{s=\pm 1}=
\frac{\partial f}{\partial \Gamma}\bigg|_{s=\pm 1}\, \lim_{t\to \pm \infty}\,\frac{\dot{\Gamma}(t)}{\dot{g}_{}(t)} = 0,
$$
where the second equality follows from Proposition~\ref{l2}.
\end{itemize}
\end{rmk}

\subsection{Attractors, repellers and invariant manifolds}

In the particular case of an attractor in the future limit system (and a repeller in the past limit system), this can be stated more topologically.  We say that a compact invariant set $A$ is 
an attractor if it is the $\omega$-limit set of a neighborhood of itself, i.e., 
there is an open set $D$ with $A \subset D$, so that $\omega (D) = A$. 
With the notation that $\psi(t,y_0)$ is the flow evolution 
for time $t$ of the initial condition $y_0$, which also 
applies to a set of initial conditions 
$\psi(t,D)=\{\psi(t,y): y\in D \}$, the $\omega-$limit set is given by:
$$
\omega (D)= \bigcap_{T>0}\,\overline{\{\psi(t,D):t>T\}},
$$
where $\overline{S}$ denotes the closure of $S$. Similarly, the $\alpha$-limit set is given by
$$
\alpha(D)= \bigcap_{T<0}\,\overline{\{\psi(t,D):t<T\}},
$$
and a set {$A$} is a repeller if it is the $\alpha$-limit set of a neighborhood of itself.
The following proposition is then readily proved for the future and past limits.
\begin{proposition}
\label{prop:attractor}
 Suppose that~\eqref{eq:ode} and $\Gamma(t)$ are chosen such that 
Proposition~\ref{l2} applies. If $A^+$ is an attractor for the future limit 
system~\eqref{eq:odea+}, then $\tilde{A}^+ =\{(x,1):x\in A^+\}\subset\{s=1\}$ is also an 
attractor for the compactified system~\eqref{eq:odeextb}--\eqref{eq:lambda_q} 
when considered in the extended phase space. If $A^-$ is a repeller for the 
past limit system~\eqref{eq:odea-}, then $\tilde{A}^- =\{(x,-1):x\in A^-\}\subset\{s=-1\}$ 
is a repeller for~\eqref{eq:odeextb}--\eqref{eq:lambda_q} when considered in 
the extended phase space.
\end{proposition}

\noindent
{\em Proof of Proposition~\ref{prop:attractor}}\\
It suffices to show the statement for an attractor as that for a repeller follows by reversing time.

From the theory of isolating blocks, see~\cite{Conley1971}, we can find a neighborhood of $A^+ $ inside $\{s=1\}$, which we denote $B$, so that the flow takes points on the boundary $\partial B$ to the interior of $B$ for $t>0$. Moreover, since this is stable under the continuous perturbation of vector fields and $\dot{s} >0 $ for $s\in(-1,1)$, we can construct a neighborhood of $\tilde{A}^+$, namely $\tilde{B}=B \times (1-\delta,1]$ for some $\delta>0$ that is forward invariant. 

Since $A$ is the maximal invariant set in $B$ and $\omega(\tilde{B})\subset B$ we can conclude that $\omega(\tilde{B}) = \tilde{A}^+$ and thus $\tilde{A}^+$ is an attractor in the full compactified system. 
\qed\\

 In the following, it will be important to distinguish between the stable set and stable manifold of a compact invariant set (similarly the unstable set and unstable manifold.)
For any compact invariant set {$A$}, we define its stable set $w^s({A})$
and unstable set $w^u({A})$ as
$$
w^s({A}) = \{y:\omega(y)\subseteq{A}\},\;\;
w^u({A}) = \{y:\alpha(y)\subseteq{A}\}.
$$
The extra condition for the stable manifold is that the decay is exponential, i.e., the stable 
manifold $W^s ({A})$ and unstable manifold $W^u ({A})$ are given by
\begin{align}
W^s ({A})=\{y \in w^s({A}) : d(\psi(y , t), {A}) \le K e^{\beta t}, \text{ for some } K>0, \beta <0  \text{ and all }t>0 \},\nonumber\\
W^u ({A})=\{y \in w^u({A}) : d(\psi(y , t), {A}) \le K e^{\beta t}, \text{ for some } K>0, \beta > 0  \text{ and all }t<0 \}.\nonumber
\end{align}
The distance function $d$ is between a point and a set, and is the greatest lower bound in the Euclidean distance between the point on the trajectory and any point in the set ${A}$:
$$
d(\psi(y,t),{A}) = \inf_{x\in{A}}\,\Vert\psi(y,t) - x\Vert.
$$
We will also need local versions of the above sets/manifolds. Based on an open superset
$N$ of a 
compact invariant set ${A}$, we define the local stable set of ${A}$:
\begin{align}
 w^s_{loc}({A})= \{y:y\in w^s({A}) \text{ and }\psi(y,t)\in N \text{ for all } t>0 \},\nonumber
\end{align}
and a local stable manifold of ${A}$:
\begin{align}
 W^s_{loc}({A})= \{y:y\in W^s({A}) \text{ and }\psi(y,t)\in N \text{ for all } t>0 \}.\nonumber
\end{align}
The local unstable set of ${A}$, denoted $w^u_{loc}({A})$, and the local unstable manifold 
of ${A}$, denoted $W^u_{loc}({A})$, are defined similarly for $t<0$.

Due to the possibility that $l^\pm_{s}=0$, we need to consider compact invariant sets with centre directions. For now, we focus on equilibria and later on, in Sec.~\ref{sec:gencis}, generalise to more 
complicated compact invariant sets.

The corresponding invariant manifolds cannot be described in terms of decay rates so straightforwardly. They are characterised as the graphs of functions over the relevant subspaces from the linearised system.
Specifically, let $E^c,E^s$ and $E^u$ be the subspaces based at $\eta$ and spanned by the sets of eigenvectors corresponding to centre (zero real-part) eigenvalues, stable (negative real-part) eigenvalues and unstable (positive real-part) eigenvalues, respectively. These subspaces are referred to as the centre ($E^c$), stable ($E^s$) and unstable ($E^u$) eigenspaces, and are invariant under the linearised system at $\eta$.
A (local) centre manifold is given as the graph of a Lipschitz function $h^c$ that is flow-invariant relative to some chosen neighborhood of the equilibrium $\eta$. More precisely,
$$
W^c_{loc}(\eta) =\{ {\rm gr}(h^c)\hspace{.05in} {\rm where} \hspace{.05in} h^c:E^c\cap N \rightarrow E^s \oplus E^u\},
$$
where $W^c_{loc}(\eta)$ is flow-invariant relative to $N$ and tangent to $E^c$ at $\eta$. 
The (local) centre-stable manifold is defined similarly, as the graph of a Lipschitz function $h^{cs}$ that is flow-invariant relative to some 
chosen neighborhood of the equilibrium $\eta$, but with a different domain space.  More precisely,
$$
W^{cs}_{loc}(\eta)=\{ {\rm gr}(h^{cs}) \hspace{.05in} {\rm where} \hspace{.05in} h^{cs}:E^c \oplus E^s \cap N \rightarrow E^u\},
$$
where $W^{cs}_{loc}(\eta)$  is flow-invariant relative to $N$ and tangent to $E^c \oplus E^s$ at $\eta$. 
In general, neither the centre nor centre-stable manifold are unique in that there may well be other functions whose graphs satisfy the tangency conditions and are invariant relative to $N$ with the corresponding domains and ranges.

\subsection{Dynamical structure from the future limit system}
\label{sec:dsfl}

Of particular interest in both the R-tipping and radial steady state problems mentioned in Sec.\ref{sec:examples} 
is the situation where a hyperbolic saddle equilibrium, denoted ${\eta}^+$, is present for the future 
limit system. In the extended phase space of the compactified system, this saddle becomes 
\begin{equation}
    \tilde{\eta}^+ = (\eta^+,1)\in\{s=1\}.\nonumber
\end{equation}

If the decay of $\dot{\Gamma}(t)$ to zero as $t\to+\infty$ is exponential or faster, then 
it is possible to construct a transformation $s=g(t)$ [e.g. transformation~\eqref{eq:gexp1}] 
so the saddle gains an exponentially stable direction and remains hyperbolic with a higher-dimensional 
stable manifold when embedded in the extended phase space of~\eqref{eq:odeextb}--\eqref{eq:lambda_q}; 
see Corollary~\ref{corol:1}. In this case, it follows  from the stable manifold theorem 
that $W^s(\tilde{\eta}^+) = w^s(\tilde{\eta}^+)$. However, this may not be true 
in general.

\begin{figure}[t]
  \begin{center}
    \includegraphics[width=14cm]{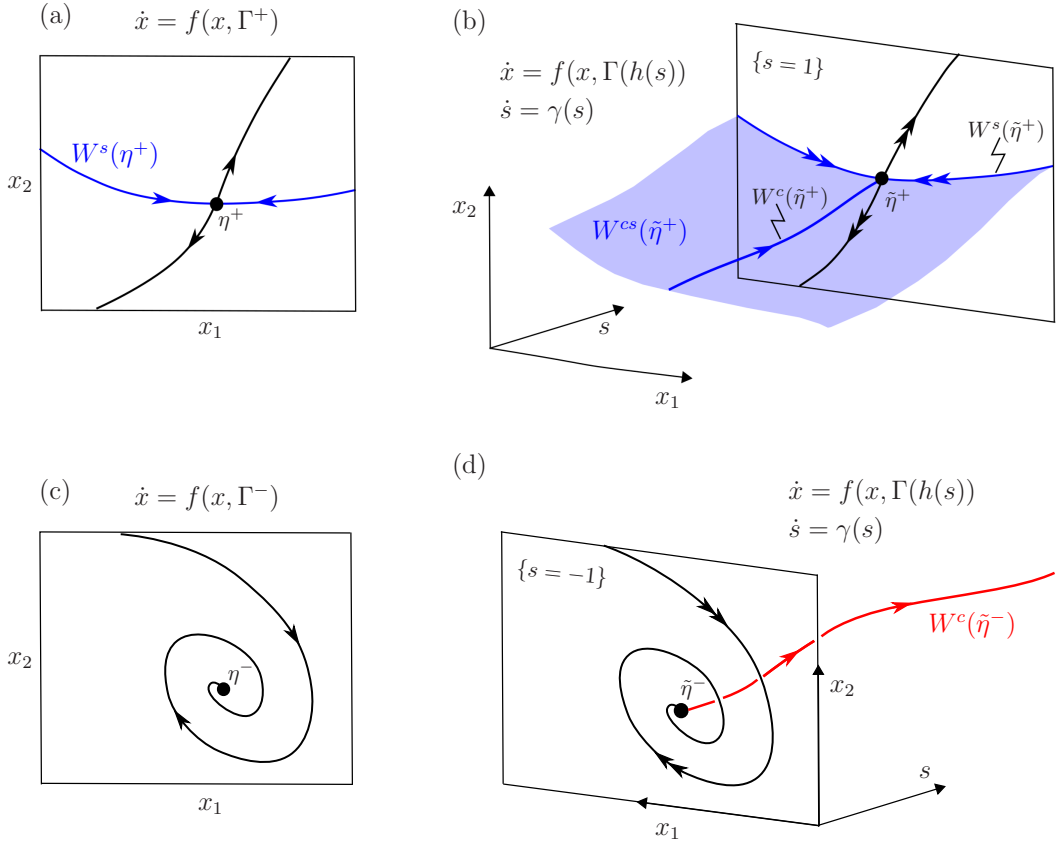}
  \end{center}
  \vspace{-3mm}
  \caption{
  (a) A hyperbolic saddle $\eta^+$ for the future limit system~\eqref{eq:odea+} 
   gains a neutrally stable centre direction and becomes (b) a non-hyperbolic saddle 
  $\tilde{\eta}^+=(\eta^+,1)$ in the extended phase space of the compactified system~\eqref{eq:odeextb}--\eqref{eq:lambda_q}.
  (c) A hyperbolic stable focus $\eta^-$ for the past limit system~\eqref{eq:odea-} 
  gains a neutrally unstable centre direction and becomes (d) a non-hyperbolic saddle-focus
  $\tilde{\eta}^-=(\eta^-,-1)$ with one-dimensional centre manifold in the extended phase 
  space of~\eqref{eq:odeextb}--\eqref{eq:lambda_q}.
  }
  \label{fig:centre_man}
\end{figure}

An interesting situation occurs when the decay of $\dot{\Gamma} (t)$ as $t\to +\infty$ is 
slower than exponential. In this case, $l_s^+ = 0$ and the saddle gains 
a neutrally stable centre direction when embedded in the extended phase space of~\eqref{eq:odeextb}--\eqref{eq:lambda_q}, meaning that 
$W^s(\tilde{\eta}^+)\neq w^s(\tilde{\eta}^+)$ because $W^s(\tilde{\eta}^+)\subset \{s=1\}$; 
see Corollary~\ref{corol:2} and Fig.~\ref{fig:centre_man}(a)--(b).
Nonetheless, we can work with a local centre-stable manifold $W^{cs}_{loc}(\tilde{\eta}^+)$ instead.

Because of the special structure inherent in the compactified system,
$s(t)$ tends monotonically to $1$,
$W^{cs}_{loc}(\tilde{\eta}^+)$ is forward invariant by construction,
and we can ensure that $W^{cs}_{loc}(\tilde{\eta}^+)\subset w^s(\tilde{\eta}^+)$. 
What is more, it turns out that the centre-stable manifold of $\tilde{\eta}^+$ actually comprises its entire local stable set, which is very useful in applications where the decay is slower 
than exponential. This follows from the following theorem. 
\begin{theorem}
\label{prop:saddler}
Consider a hyperbolic saddle equilibrium ${\eta}^+$ for the future 
limit system~\eqref{eq:odea+} that becomes a non-hyperbolic saddle 
$\tilde{\eta}^+=(\eta^+,1)$ in the compactified system~\eqref{eq:odeextb}--\eqref{eq:lambda_q}.
Then, 
$W^{cs}_{loc}(\tilde{\eta}^+)$ is unique relative to 
 some chosen neighborhood $N$ of $\tilde{\eta}^+$, and

$$
W^{cs}_{loc}(\tilde{\eta}^+) = w^s_{loc}(\tilde{\eta}^+),
$$
in the extended phase space of the compactified system~\eqref{eq:odeextb}--\eqref{eq:lambda_q}.
\end{theorem}
%
%

\noindent
{\em Proof of Theorem~\ref{prop:saddler}}\\
\noindent
Without loss of generality, we set $\eta^+=0$ and introduce $\sigma =s-1$. We can choose coordinates $x=(y,z)$ so that, in a neighborhood of $(x,\sigma)=(0,0)$, (\ref{eq:odeext})-(\ref{eq:lambda_q+}) becomes
\begin{align}
  \label{eq:odelimstable+}
  \dot{y}&= M^+y+f_1(y,z,\sigma),\\
  \label{eq:odelimunst+}
  \dot{z}&=M^-z + f_2(y,z,\sigma),\\
  \label{eq:gamnewcoords+}
  \dot{\sigma}&=  \beta(\sigma),
 \end{align}
where $\beta (\sigma)=\gamma(\sigma +1)$, and $\beta (0)= \beta'(0)=0$. The spectra of $M^+$ and $M^-$ are in the right and left half planes respectively. Set 
$$
N=\{(y,z,\sigma): \Vert y\Vert \le \epsilon, \Vert z \Vert<\epsilon, -\delta < \sigma \le 0\}.
$$ 
Since the terms $f_i$ are higher order and $\dot{\sigma}>0$, we can choose $\epsilon$ and $\delta$ so that the vector field of (\ref{eq:odelimstable+})-(\ref{eq:gamnewcoords+}) points into $N$ on $\Vert z\Vert=\epsilon$ and $\sigma=-\delta$ while it points out of $N$ on $\Vert y\Vert=\epsilon$, and the last face $\sigma =0$ ($s=1$) is invariant.

The variables $z$ and $\sigma$ represent the stable and centre directions respectively and so we put them together and introduce $\zeta=(z,\sigma)$.
We also denote $\xi =(y,\zeta)$ and 
$\psi (\xi ,t)=(y(t),\zeta(t))$.
Note that if $\xi \in w^s_{loc}(\tilde{\eta}^+)$ then $\psi (\xi ,t) \in N$ for all $t>0$ by construction
of $w^s_{loc}(\tilde{\eta}^+)$.
If $\xi \in W^{cs}_{loc}(\tilde{\eta}^+)$ then $\psi (\xi ,t) \in N$ for all $t>0$ by construction  of $W^{cs}_{loc}(\tilde{\eta}^+)$ and the direction of the vector field on the faces of $N$. 
To see this, note that $W^{cs}_{loc}(\tilde{\eta}^+)$ cannot intersect $\Vert y\Vert=\epsilon$
because it is invariant,  and the graph of a Lipschitz function $y=h^{cs}(\zeta)$ on which $\dot{\sigma}>0$.

Suppose that $\xi_1=(y_1,\zeta_1)$ 
satisfies $\psi (\xi_1,t) \in N$ for all $t\ge 0$. 
Let $\xi_2=(y_2,\zeta_1)$ be the point on $W^{cs}_{loc}$ with the same $\zeta$ coordinate as $\xi_1$, which will also stay in $N$ for $t>0$.
Under these conditions, the following lemma describes the fate of the hypothesized trajectory from $\xi_1$.

\begin{lemma}
\label{lemma:expgrowth}
There are coordinates for $y$ and $\zeta$ so that the following holds for all $t\ge0$:
$$ \Vert y_1(t)-y_2(t)\Vert\ge \Vert\zeta_1(t)-\zeta_2(t)\Vert.
$$
Furthermore, there is an $\eta >0$ so that
$$
\Vert y_1(t)-y_2(t)\Vert\ge \Vert y_1(0)-y_2(0)\Vert \exp{(\eta t)}
$$ for all $t \ge 0$.
\end{lemma}
If $y_1(0) \ne y_2(0)$ then Lemma~\ref{lemma:expgrowth} implies that $\xi_1 (t)$ must leave $N$ for some $t>0$, which is a contradiction.
%
%
In other words, if we suppose $\xi_1\in W^{cs}_{loc}(\tilde{\eta})$ 
then the contradiction gives a unique function $y=h^{cs}(\zeta)$ whose graph satisfies the tangency at $\tilde{\eta}$ and flow-invariance relative to $N$
conditions. 
If we suppose $\xi_1\in w^{s}_{loc}(\tilde{\eta})$ 
then the contradiction gives $\xi_1 \in W^{cs}_{loc}$ and  $W^{cs}_{loc}(\tilde{\eta}) = w^{s}_{loc}(\tilde{\eta})$.
Once we prove Lemma~\ref{lemma:expgrowth} the  proof of Theorem~\ref{prop:saddler} is complete.
\vspace{3mm}

\noindent
{\em Proof of Lemma~\ref{lemma:expgrowth}}\\
\noindent
It will be convenient to combine (\ref{eq:odelimunst+}) with (\ref{eq:gamnewcoords+}), and write
\begin{equation}
\label{eq:combstagam}
 \dot{\zeta}=M^0
\zeta +\tilde{f}_2(y,\zeta) .  
\end{equation}
The spectrum of $M^0$, denoted $\sigma (M^0)$, contains one neutral and otherwise stable eigenvalues, while that of $M^+$ contains unstable eigenvalues. Therefore we can find $\alpha $ and $\beta$ so that, with suitable coordinates, $\langle M^+y,y \rangle \ge \beta \langle y,y\rangle$ and $\langle M^0 \zeta, \zeta \rangle \le \alpha \langle \zeta,\zeta \rangle$, where $0<\alpha<\beta, {\rm and  } \langle \hspace{1mm},\rangle$ denotes the inner product in those coordinates. 

Now set $\Delta y =y_1-y_2$ and $\Delta \zeta =\zeta_1-\zeta _2$. The functions $f_1$ and $\tilde{f}_2$ being higher order means that $Df_1(0,0,0)=Df_2(0,0,0)=0$. From the continuity of the  first derivatives of these nonlinear terms, given $\rho >0$, we can find $\epsilon >0$ and $\delta >0$ so that
$$
\left\Vert f_1(y_1,\zeta_1)-f_1(y_2,\zeta_2) \right\Vert \le \rho \left( \Vert\Delta y\Vert + \Vert\Delta \zeta \Vert \right).
$$
A similar estimate holds for $\tilde{f}_2$. The following inequality can then be shown to hold:

\begin{align}
\label{eq:cone1}
\frac{d}{dt} \Vert\Delta y \Vert ^2 & \ge \beta \Vert\Delta y \Vert^2 -\rho(\Vert\Delta y \Vert + \Vert\Delta \zeta \Vert) \Vert\Delta y\Vert,\\
\label{eq:cone2}
\frac{d}{dt} \Vert\Delta \zeta \Vert ^2 &  \le \alpha \Vert\Delta \zeta \Vert^2 +\rho ( \Vert\Delta y\Vert + \Vert\Delta \zeta \Vert )\Vert\Delta \zeta\Vert.
\end{align}
Subtracting (\ref{eq:cone2}) from (\ref{eq:cone1}), we obtain:
$$
\frac{d}{dt} (\Vert\Delta y \Vert ^2 - \Vert\Delta \zeta \Vert ^2) \ge \beta \Vert\Delta y \Vert ^2 -\alpha \Vert\Delta \zeta \Vert^2  - \rho (\Vert\Delta y \Vert + \Vert\Delta \zeta \Vert)^2.
$$
Initially, $\Vert\Delta y\Vert > \Vert\Delta \zeta \Vert$ and if they ever became equal, we can conclude from this inequality that 
$$
\frac{d}{dt} (\Vert\Delta y \Vert ^2 - \Vert\Delta \zeta \Vert ^2) \ge (\beta -\alpha - 4 \rho)\Vert\Delta y \Vert ^2.
$$
and so, if $\rho$ is sufficiently small,  $\frac{d}{dt} (\Vert\Delta y \Vert ^2 - \Vert\Delta \zeta \Vert ^2) \ge 0$,
and hence $\Vert\Delta y (t)\Vert > \Vert\Delta \zeta (t)\Vert$ for all $t>0$ so long as both trajectories stay in $N$, which they do by hypothesis.

It is now straightforward to obtain the second inequality of the lemma from (\ref{eq:cone1}) using $\Vert\Delta y (t)\Vert > \Vert\Delta \zeta (t)\Vert$ and setting $\eta = \beta- 2 \rho$.
\qed

\subsection{Dynamical structure from the past limit system}
\label{sec:dspl}

Analogous results can be proved for the past limit.
Of particular interest in the pullback attractor and R-tipping problems mentioned
in Sec.\ref{sec:examples}  is the situation where a hyperbolic sink, denoted ${\eta}^-$, 
is present for the past limit system. In the extended phase space of the compactified 
system, this equilibrium becomes 
\begin{equation}
    \tilde{\eta}^- = (\eta^-,-1) \in\{s=-1\}.\nonumber
\end{equation}
 If the decay of $\dot{\Gamma}(t)$ to zero as $t\to-\infty$ 
is exponential or faster, then it is possible to construct a transformation $s=g(t)$ [e.g. transformation~\eqref{eq:gexp1}] so the sink gains an exponentially unstable direction and becomes 
a hyperbolic saddle with a one-dimensional unstable manifold when embedded in the extended phase space of~\eqref{eq:odeextb}--\eqref{eq:lambda_q}; see Corollary~\ref{corol:1}.
In this case, it follows  from the stable manifold theorem
that $W^u(\tilde{\eta}^-) = w^u(\tilde{\eta}^-)$. However, this may not be true 
in general.

An interesting situation occurs when the decay of $\dot{\Gamma} (t)$ as $t\to -\infty$ is 
slower than exponential. In this case the sink gains 
a neutrally unstable centre direction when embedded in the extended phase space of~\eqref{eq:odeextb}--\eqref{eq:lambda_q}, meaning that 
$W^u(\tilde{\eta}^-)\neq w^u(\tilde{\eta}^-)$ because $W^u(\tilde{\eta}^-)=\emptyset$; 
see Corollary~\ref{corol:2} and Fig.~\ref{fig:centre_man}(c)--(d).
Nonetheless, we can work with a local centre manifold $W^{c}_{loc}(\tilde{\eta}^-)$.

Because of the special structure inherent in the compactified system,
$s(t)$ tends monotonically to $-1$ as $t\to - \infty$,
$W^{c}_{loc}(\tilde{\eta}^-)$ is backward invariant by construction,
and we can ensure that $W^{c}_{loc}(\tilde{\eta}^-)\subset w^u(\tilde{\eta}^-)$. 
What is more, it turns out that the centre manifold of $\tilde{\eta}^-$ actually comprises its entire local unstable set,
which is very useful in applications where the decay is slower 
than exponential. This follows from the following theorem. 
\begin{theorem}
\label{prop:sinkr}
Consider a hyperbolic sink ${\eta}^-$ for the past limit system~\eqref{eq:odea-}
that becomes a non-hyperbolic saddle $\tilde{\eta}^-=(\eta^-,-1)$ in the compactified system~\eqref{eq:odeextb}--\eqref{eq:lambda_q}. Then, 
$W^{c}_{loc}(\tilde{\eta}^-)$ is unique relative to
some chosen neighborhood $N$ of $\tilde{\eta}^-$, and

$$
W^{c}_{loc}(\tilde{\eta}^-) = w^u_{loc}(\tilde{\eta}^-),
$$
in the extended phase space of the compactified system~\eqref{eq:odeextb}--\eqref{eq:lambda_q}.
\end{theorem}
\noindent
{\em Proof of Theorem~\ref{prop:sinkr}}\\
This follows immediately from Theorem~\ref{prop:saddler} by reversing time and seeing that $W^{cu}_{loc}(\tilde{\eta}^-) = w^u_{loc}(\tilde{\eta}^-)$. But, in this case, there is no unstable direction and so $W^{cu}_{loc}=W^{c}_{loc}$.
\qed

\subsection{Invariant manifolds and applications}

The key point of Sec.~\ref{sec:compactdyns} is that the unique invariant manifolds extrapolated 
from the future limit system, from the past limit system, or from both, 
can be used to analyse the nonautonomous system~\eqref{eq:ode} 
with exponential or non-exponential
decaying $\Gamma(t)$. Specifically, the benefits of compactification when studying
nonlinear dynamics are the following:
 \begin{itemize}
\item[(i)] In contrast to the original nonautonomous system~\eqref{eq:ode}, the 
compactified system~\eqref{eq:odeextb}--\eqref{eq:lambda_q} is autonomous and may 
contain equilibria and more complicated compact invariant sets, e.g. limit cycles, invariant tori or strange attractors.
\item[(ii)] 
In the extended phase space of~\eqref{eq:odeextb}--\eqref{eq:lambda_q}, the temporal or spatial variation of the  nonautonomous term $\Gamma(t)$ 
becomes in a certain sense `encoded' in  the geometric shape of: 
\begin{itemize}
    \item[$\bullet$] The stable or centre-stable invariant manifold of a saddle $\tilde{\eta}^{+}$.
    \item[$\bullet$] The unstable or centre invariant manifold of a saddle $\tilde{\eta}^{-}$.
    \item[$\bullet$] The basin of attraction of an attractor  $\tilde{A}^+$.
\end{itemize}
\item[(iii)] 
In the extended phase space of~\eqref{eq:odeextb}--\eqref{eq:lambda_q}:
\begin{itemize}
    \item[$\bullet$] Orbits contained either in the basin of attraction of an
    attractor $\tilde{A}^+$, or in the stable or centre-stable invariant  manifold of a saddle $\tilde{\eta}^{+}$, correspond to solutions of the nonautonomous system~\eqref{eq:ode} that remain bounded as $t\to +\infty$.
    \item[$\bullet$] Orbits contained in the unstable or centre invariant manifold of a saddle $\tilde{\eta}^{-}$  correspond to solutions of the nonautonomous system~\eqref{eq:ode} that remain bounded as $t\to -\infty$.
    \item[$\bullet$] Connecting heteroclinic orbits from a saddle $\tilde{\eta}^{-}$ to either a saddle $\tilde{\eta}^{+}$ or an attractor $\tilde{A}^+$ correspond to solutions of the nonautonomous system~\eqref{eq:ode} that  remain bounded as $t$ tends to
    positive and negative infinity.
\end{itemize}
\item[(iv)] 
The compactification enables numerical detection and parameter continuation of 
solutions to the nonautonomous system~\eqref{eq:ode} that  remain bounded as $t\to\pm\infty$. 
Since such solutions become connecting heteroclinic orbits in the compactified 
system~\eqref{eq:odeextb}--\eqref{eq:lambda_q}, they can be detected and continued numerically, 
for example using the implementation of the Lin's method~\cite{Lin1990} in the numerical 
continuation software AUTO~\cite{auto,Oldeman2003,Krauskopf2008,Xie2019}.
\end{itemize}

These benefits are best illustrated using the three motivating examples from Sec.~\ref{sec:examples}.

\subsubsection{Pullback attractors}
A local pullback attractor for a nonautonomous system that is asymptotically constant 
as $t\to\ - \infty$ will emanate from a normally hyperbolic attractor $A^-$ for the past 
limit system. 
For example, consider a sink $\eta^-$ for the past limit system, which becomes a
saddle $\tilde{\eta}^{-}=(\eta^-,-1)$ when embedded in the extended phase space of the 
compactified system~\eqref{eq:odeextb}--\eqref{eq:lambda_q}. Then, the local pullback 
point attractor of~\eqref{eq:ode} associated with $\eta^-$ is transformed into the unstable 
or centre invariant manifold of the saddle $\tilde{\eta}^{-}$ in the extended phase space of~\eqref{eq:odeextb}--\eqref{eq:lambda_q}.

\subsubsection{Rate-induced tipping}
In the R-tipping problems, the objective is to find thresholds for the tipping. We claim these are naturally anchored in the future limit system by a saddle or repeller that is normally hyperbolic. We call such a state
for the future limit system an {\em R-tipping edge state}; see~\cite{Wieczorek2019} where this terminology is fully explained. 
The challenge is to pin down, in the extended phase space of~\eqref{eq:odeextb}--\eqref{eq:lambda_q}, what we call {\em R-tipping thresholds}~\cite{Wieczorek2019} separating solutions of~\eqref{eq:ode} that R-tip from those that do not. The key point is that these can be naturally characterised as the $n$-dimensional stable or centre-stable manifolds of an appropriate saddle edge state. 
What is more, the critical rates of change of $\Gamma(t)$ at which the system transitions 
to a different state correspond to heteroclinic connections from a saddle $\tilde{\eta}^{-}$ to 
the saddle edge state. 
Thus, the compactification simplifies rigorous analysis, and enables numerical detection and parameter continuation of genuine nonautonomous R-tipping bifurcations~\cite{Xie2019}.

\subsubsection{Radial steady states}
For the construction of radial steady states as in (\ref{eq:sphwaves}), the invariant manifolds constructed above will provide the sets of solutions satisfying the boundary conditions. We would be interested in saddles for the future 
limit system. This is because,  in the extended phase space of~\eqref{eq:odeextb}--\eqref{eq:lambda_q}, the (family of) solutions 
to the radial steady state problem~\eqref{eq:ode} must decay as $r\rightarrow +\infty$. This would then use a right-sided compactification.  A separate construction is needed for the boundary condition at $r=0$. When this has been carried out, the problem becomes one of finding heteroclinic connections given by transverse intersections of the unstable or centre manifold of a saddle $\tilde{\eta}^{-}$ and the stable or centre-stable manifold of a saddle $\tilde{\eta}^{+}$, see~\cite{Jones1986}.

\subsubsection{Stability of nonlinear waves}
The Evans Function is defined using sets of solutions of (\ref{eq:twev})-(\ref{eq:twev2}) which decay at $\pm \infty$, denoted $Y_\pm(z)$. Using exterior algebra, the Evans Function is defined as a generalised angle between these subspaces when compared at a fixed $z \in (-\infty, +\infty)$. These subspaces are realised in the compactified system as (at $-\infty$) a centre-unstable manifold of the zero solution for the left-sided compactification and, analogously, (at $+\infty$) a centre-stable manifold of the zero solution for the right-sided compactification; see~\cite{Alexander1990}. 

In most circumstances, the decay of the underlying wave to its end-states $U_\pm$
is exponential. Nevertheless, the theory presented here would afford a construction of the subspaces $Y_\pm (z)$, and hence the Evans Function, when the decay were much weaker than exponential.

\subsection{Generalisations to more complicated compact invariant sets}
\label{sec:gencis}
We anticipate that it will be possible to extend the results {\ for equilibria} from Secs.~\ref{sec:dsfl} and~\ref{sec:dspl} to cover more complicated compact invariant sets in the limit systems and their invariant manifolds in the 
extended phase space of the compactified system. It is worth distinguishing two separate cases here.

If the compact invariant set $A$ is a submanifold of $\{s=1\}$, then the relevant invariant manifold theorems can be inferred from the results in Fenichel \cite{Fenichel1971}. Theorems~\ref{prop:attractor} and \ref{prop:saddler} will follow in a straightforward manner with the use of sub-bundles of the relevant tangent spaces as opposed to subspaces. 

Arnol'd et al. \cite[Sec.8.2]{Arnold2013} define a strange attractor as an attractor that is not the union of smooth submanifolds. Similarly, we might define a strange compact invariant set as one that cannot be cast as the union of submanifolds. The idea from the previous paragraph should extend to the case of a union of submanifolds, and so all non-strange compact invariant sets will be covered. This will obviously include important invariant sets such as periodic orbits and tori. 

Invariant manifold theorems do exist for general compact invariant sets, including strange ones, see the book by Robinson \cite{Robinson1999} for instance. Nevertheless, it is not apparent how to generalise the arguments in the proofs of Theorems~\ref{prop:saddler} and \ref{prop:sinkr} as these proofs require a nicer local structure than might hold near a more complicated compact invariant set or strange attractor. We formulate the general result as a conjecture.

\begin{conj}
\label{con:highdim}
The statements for equilibria in Theorems~\ref{prop:saddler} 
and~\ref{prop:sinkr} also hold for general normally hyperbolic compact invariant sets 
{$A^\pm$} of the limit systems. 
\end{conj}

\section{Examples of Compactification}
\label{sec:AppA}

In this section we  construct actual examples of coordinate 
transformations $g(t)$, given by~\eqref{eq:g},~\eqref{eq:g+} and~\eqref{eq:g-} 
and depicted in Fig.~\ref{fig:g}(a)--(c), that are useful in practice. 
This raises the question of the `rate' of compactification relative to the 
rate of asymptotic decay of the nonautonomous term. To address this question, we introduce 
in Secs.~\ref{sec:ec} and~\ref{sec:ac}
{compactification parameters} and work with {\em parametrised compactifications}.
The additional freedom acquired in a parametrisation of $g(t)$ allows us to:
\begin{itemize}
\item[(i)] Fulfill the compactification conditions~\eqref{eq:dLdg}--\eqref{eq:d2gdt2} by a suitable choice of the compactification parameter(s).
\item[(ii)] Control the direction of the additional Lyapunov vector 
transverse to invariant subspaces $\{ s=\pm 1 \}$ to facilitate analysis.
\end{itemize}
Moreover, we consider different types of asymptotic decay of $g(t)$ giving rise to 
different behaviour of $s(t)$ near $\{s=\pm 1\}$.

\subsection{$\Gamma$-Compactification}
\label{sec:Lc}

In some cases, a compactification transformation can be constructed in 
terms of (one component of) the nonautonomous term $\Gamma(t)\in\R^d$ itself; see~\cite{Ashwin2012,Perryman2014,Perryman2015,Ritchie2016,Alkhayuon2018} 
for examples. 
Let $\dot{\Gamma}(t)=(\dot{\Gamma}_1(t),\ldots,\dot{\Gamma}_d(t))$
be ordered such that  $\dot{\Gamma}_1(t)$ decays not faster than the other components
$$
\lim_{t\to\pm\infty}\frac{\Gamma_i(t)}{\Gamma_1(t)} = L^\pm\in\R\;\;\mbox{for}\;\;i=2,\dots,d.
$$

If $\Gamma_1(t)$ satisfies Assumption~\ref{asmp:g} with 
future limit $\Gamma_1^+$ and past limit $\Gamma_1^-$, we can construct a  coordinate transformation~\eqref{eq:g} in terms of $\Gamma_1(t)$ as follows
\begin{equation}
\label{Gtrans}
g(t)=\frac{2\Gamma_1(t) - \Gamma_1^+ - \Gamma_1^-}{\Gamma_1^+ -
  \Gamma_1^-}.
\end{equation}
One can verify that such $g(t)$ satisfies
the first transformation condition~(\ref{eq:dLdg}):
\begin{equation}
\label{eq:Lc2}
\lim_{t\to \pm\infty}\frac{\dot{\Gamma}_i(t)}{\dot{g}(t)}=
\left\{\begin{array}{rcl}
      (\Gamma_1^+ - \Gamma_1^-)/2\ne 0  &\mbox{for} & i= 1,\nonumber\\
     (\Gamma_1^+ - \Gamma_1^-)/2\ne 0\;\;\mbox{or}\;\;0 &\mbox{for} & i=2,\ldots,d,
\end{array}
\right.
\end{equation}
and satisfies the second transformation condition~(\ref{eq:d2gdt2}):
$$
\lim_{t\to \pm\infty}\frac{\ddot{g}(t)}{\dot{g}(t)}=
\lim_{t\to \pm\infty}\frac{\ddot{\Gamma}_1(t)}{\dot{\Gamma}_1(t)}\;\;\mbox{exists},
$$
if $\Gamma_1(t)$ satisfies~(\ref{eq:d2gdt2}).
Similarly, one can construct a coordinate transformation~\eqref{eq:g+}:
$$
g(t)=\frac{\Gamma_1(t) - \Gamma_1(0)}{\Gamma_1^+ - \Gamma_1(0)},
$$ 
if $\Gamma_1(t)$ satisfies Assumption~\ref{asmp:g+} with future limit $\Gamma_1^+$, 
or a coordinate transformation~\eqref{eq:g-}:
$$
g(t)= - \frac{\Gamma_1(t) - \Gamma_1(0)}{\Gamma_1^- - \Gamma_1(0)},
$$ 
if $\Gamma_1(t)$ satisfies Assumption~\ref{asmp:g-} with past limit $\Gamma_1^-$.

For example, one can use~\cite{Ashwin2012,Perryman2014}:
$$g(t)=\tanh(t)\quad\mbox{if}\quad \Gamma_1(t) = \tanh(t).$$ 
This approach can be extended to certain non-monotone inputs that can be expressed in terms of a monotone growth that satisfies Assumption~\ref{asmp:g}, such as $$\Gamma_1(t) = \mbox{sech}(t) =\sqrt{1-\tanh(t)^2}.$$

Constructing $g(t)$ in terms of $\Gamma_1(t)$ is not always possible, for example it usually
does not work for non-monotone $\Gamma_1(t)$. When possible, it has certain 
limitations. For example, there may be no algebraic formula for the inverse $h(s)$ 
required for the augmented component~\eqref{eq:odeext2b} of the vector field. 
Moreover, to facilitate analysis, one may wish 
the additional Lyapunov vectors to be normal 
to invariant subspaces $\{s=\pm 1\}$, which does not normally occur for the
$\Gamma$-compactification because the resulting Jacobian~\eqref{eq:jac} is not block-diagonal. Therefore, more universal
compactification transformations will be required in general.

\subsection{Exponential Compactification}
\label{sec:ec}

Guided by Assumption~\ref{asmp:g} and Proposition~\ref{l2}, we construct 
an example of a {\em parametrised coordinate transformation}~(\ref{eq:g}) with 
exponential asymptotic decay as proposed in~\cite[Eq.(3.2)]{Alexander1990}:
\begin{align}
\label{eq:gexp1}
s = g_{(\alpha)}(t) =\tanh\left(\frac{\alpha t}{2}\right),\;\; 
t = g^{-1}_{(\alpha)}(s):= h_{(\alpha)}(s) = \frac{1}{\alpha}\ln\frac{1+s}{1-s},
\end{align}
where the {\em compactification parameter}  $\alpha > 0$ parametrises the {\em rate} of exponential decay of $\dot{g}_{(\alpha)}(t)$.
One can  verify that this transformation satisfies Assumption~\ref{asmp:g}  and
the second transformation condition~\eqref{eq:d2gdt2}:
\begin{equation}
\label{eq:tc2exp}
\lim_{t\to \pm\infty}\,\frac{\ddot{g}_{(\alpha)}(t)}{\dot{g}_{(\alpha)}(t)} =
-\alpha \lim_{t\to \pm\infty} \tanh\left(\frac{\alpha t}{2}\right) =\mp \alpha.
\end{equation}
The augmented component $\gamma_{(\alpha)}(s)$ of the vector field  is obtained
by computing and then inverting 
$h_{(\alpha)}'(s)$,
so that the compactified system \eqref{eq:odeextb}--\eqref{eq:odeext2b} becomes
\begin{eqnarray}
\label{eq:dsdtexp}
\left.
\begin{split}
\dot{x} &= f(x,\Gamma(h_{(\alpha)}(s))),\\
\dot{s} &=  \frac{\alpha}{2} \left(1 - s^2\right).
\end{split}
\right\}
\end{eqnarray}
We now examine properties of~\eqref{eq:dsdtexp} with dependence on $\alpha$ and the asymptotic decay of $\dot{\Gamma}(t)$.
\begin{corol}
\label{corol:1} {\em (Compactification for $\dot{\Gamma}(t)$ with exponential or faster decay.)}
Consider a nonautonomous system~(\ref{eq:ode}) with $C^1$-smooth $f$ and $\Gamma$, and 
 bi-asymptotically constant $\Gamma(t)$. Moreover, suppose there is a $\rho>0$ such that
\begin{align}
\label{eq:expdec}
&\lim_{t\to \pm\infty}\frac{\dot{\Gamma}(t)}{e^{\mp\rho t}}\;\;\mbox{exist}.
\end{align}
Then, given the coordinate transformation~\eqref{eq:gexp1} with any 
$\alpha\in(0,\rho]$, the ensuing autonomous 
compactified system~\eqref{eq:dsdtexp} is $C^1$-smooth on the extended phase space 
$U\times[-1,1]$. \\
The distance between any trajectory of~\eqref{eq:dsdtexp} and $\{s= \pm1\}$ decays
exponentially at a rate $\alpha$ as $t\to\pm\infty$.
A regular compact invariant set of the future limit system gains one negative Lyapunov exponent
$l_{s,(\alpha)}^+ = -\alpha<0$, and a regular compact invariant set of the past limit system gains one positive Lyapunov exponent
$l_{s,(\alpha)}^- = \alpha>0$, when embedded in the extended phase space of~\eqref{eq:dsdtexp}.
\end{corol}
In other words, exponentially or faster decaying $\dot{\Gamma}(t)$ is sufficient for~\eqref{eq:dsdtexp}
to be $C^1$-smooth and for regular compact invariant sets of the limit systems to gain a non-zero Lyapunov exponent  
when embedded in the extended phase space of~\eqref{eq:dsdtexp}. This gives rise to pure stable and unstable manifolds in the extended phase space.
\vspace{3mm}

\noindent
{\em Proof of Corollary~\ref{corol:1}}\\
\noindent
Given~\eqref{eq:expdec}, it follows from the asymptotic properties of 
$\dot{g}_{(\alpha)}(t)$:
\begin{align}
\dot{g}_{(\alpha)}(t)&\sim 2\alpha \, e^{\mp\alpha t}\;\;\mbox{as}\;\;
t\to\pm\infty,\nonumber
\end{align}
and from the algebraic limit theorem
\begin{align}
\lim_{t\to \pm\infty}\frac{\dot{\Gamma}(t)}{\dot{g}_{(\alpha)}(t)}&= 
\lim_{t\to \pm\infty}\frac{\dot{\Gamma}(t)}{e^{\mp\rho t}} \frac{e^{\mp\rho t}}{\dot{g}_{(\alpha)}(t)}=
\lim_{t\to \pm\infty}\frac{\dot{\Gamma}(t)}{e^{\mp\rho t}}
\lim_{t\to \pm\infty}\frac{e^{\mp\rho t}}{\dot{g}_{(\alpha)}(t)}=
 \frac{1}{2\alpha}\lim_{t\to \pm\infty}\frac{\dot{\Gamma}(t)}{e^{\mp\rho t}} 
\lim_{t\to \pm\infty}e^{\mp(\rho-\alpha) t},\nonumber
\end{align}
 that the first transformation condition~\eqref{eq:dLdg} is satisfied if $0<\alpha\le\rho$.
It follows from~\eqref{eq:tc2exp} that the second transformation
condition~\eqref{eq:d2gdt2} is satisfied for $\alpha > 0$.
It then follows from Proposition~\ref{l2} that the compactified system~\eqref{eq:dsdtexp} 
is $C^1$-smooth on the extended phase space $U\times[-1,1]$ for any $0<\alpha\le\rho$. 

The Hausdorff semi-distance $d(p(t),S^\pm)$ between a point  $p(t)=(x(t),s(t))$ on a trajectory and invariant subspace $S^\pm=\{(x,s):s=\pm 1\}$ in the extended phase space 
is the distance between $s(t)$ and $\pm 1$. Thus, exponential decay of $d(p(t),S^\pm)$ follows from exponential approach of $s(t)$ 
towards $\pm1$ as $t\to\pm\infty$:
$$
d(p(t),S^\pm) = \inf_{y\in S^\pm}\,\Vert p(t) - y\Vert  =  1\mp s(t) = \frac{2}{e^{\pm\alpha t}+1}.
$$
It folows from Eq.~\eqref{eq:sLyap} that the additional Lyapunov exponent $l_{s,(\alpha)}^{\pm}$
is given by Eq.~\eqref{eq:tc2exp}.
\qed\\

Similarly, one can discuss one-sided subcases of the two-sided compactification above,
using an example of parametrised coordinate transformation~\eqref{eq:g+} with exponential decay in the form
\begin{equation}
s = g_{(\alpha)}(t)= 1- e^{-\alpha t},\;\; t = h_{(\alpha)}(s) = - \frac{1}{\alpha}\ln(1-s),
\end{equation}
or example of parametrised coordinate transformation~\eqref{eq:g-} with exponential decay in the form
\begin{equation}
s = g_{(\alpha)}(t)= -1 + e^{\alpha t},\;\; t = h_{(\alpha)}(s) =  \frac{1}{\alpha}\ln(1+s).
\end{equation}

\subsection{Algebraic Compactification}
\label{sec:ac}

Guided by Assumption~\ref{asmp:g} and Proposition~\ref{l2}, we construct an inverse of the {\em parametrised compactification transformation}~(\ref{eq:g}) with algebraic asymptotic decay in the form akin to stereographic projection
\begin{align}
\label{eq:halg1}
t = g^{-1}_{(\alpha)}(s) := h_{(\alpha)}(s)=\frac{s}{(1-s^2)^\frac{1}{\alpha}},
\end{align}
where the {\em compactification parameter}  $\alpha >0$ parametrises the {\em order} of algebraic decay of $\dot{g}_{(\alpha)}(t)$.
There is no formula for the corresponding compactification transformation $s=g_{(\alpha)}(t)$ in the general case $\alpha > 0$. However,
there is one in the special case $\alpha =1$:
\begin{align}
\label{eq:galg1}
s = g(t) = \frac{-1 + \sqrt{1 + 4t^2}}{2t},\;\;t = h(s) = \frac{s}{1-s^2}.
\end{align}
One can  verify that transformation~\eqref{eq:halg1} satisfies Assumption~\ref{asmp:g}  and
the second transformation condition~\eqref{eq:d2gdt2}:
\begin{equation}
\label{eq:tc2alg}
\lim_{s\to \pm 1}\, \frac{h_{(\alpha)}''(s)}{\left(h_{(\alpha)}'(s)\right)^2} =
\lim_{s\to \pm 1}\,\frac{2(2-\alpha)s^3 +6\alpha s }{\left[\alpha(1-s^2) + 2s^2\right]^2}\,(1-s^2)^\frac{1}{\alpha}
= 0.
\end{equation}
The augmented component $\gamma_{(\alpha)}(s)$ of the vector field  is obtained
by computing and then inverting $h_{(\alpha)}'(s)$, so that
the compactified system \eqref{eq:odeextb}--\eqref{eq:odeext2b} becomes
\begin{eqnarray}
\label{eq:dsdtalg}
\left.
\begin{split}
\dot{x} &= f(x,\Gamma(h_{(\alpha)}(s))),\\
\dot{s} &=  \frac{\alpha(1-s^2)^{1+\frac{1}{\alpha}}}{\alpha(1-s^2) + 2s^2} .
\end{split}
\right\}
\end{eqnarray}
We now examine properties of~\eqref{eq:dsdtalg} with dependence on $\alpha$ and 
the asymptotic decay of $\dot{\Gamma}(t)$.
\begin{corol}
\label{corol:2} {\em (Compactification for $\dot{\Gamma}(t)$ with algebraic or faster decay.)}
Consider a nonautonomous system~(\ref{eq:ode}) with $C^1$-smooth $f$ and $\Gamma$, and bi-asymptotically constant $\Gamma(t)$. Moreover, suppose there is an $m > 1$ such that
\begin{align}
\label{eq:algdec}
&\lim_{t\to \pm\infty}\frac{\dot{\Gamma}(t)}{t^{-m}}\;\;\mbox{exist}.
\end{align}
Then, given the coordinate transformation~\eqref{eq:gexp1} with any $\alpha\in(0,m-1]$, the ensuing autonomous 
compactified system~\eqref{eq:dsdtalg} is $C^1$-smooth on 
the extended phase space $U\times[-1,1]$.\\
The distance between any trajectory of~\eqref{eq:dsdtalg} and $\{s= \pm1\}$ decays
algebraically as $t\to\pm\infty$.
A regular compact invariant set of a limit system gains one zero Lyapunov exponent
$l_{s,(\alpha)}^\pm = 0$ when embedded in the extended phase space of~\eqref{eq:dsdtalg}.
\end{corol}
In other words, algebraically or faster decaying $\dot{\Gamma}(t)$ is sufficient for~\eqref{eq:dsdtalg} to be $C^1$-smooth, but regular compact invariant sets for the 
limit systems gain one zero Lyapunov exponent 
when embedded in the 
extended phase space of~\eqref{eq:dsdtalg}. This gives rise to centre-stable and centre manifolds in the extended phase space.
\vspace{3mm}

\noindent{\em Proof of Corollary~\ref{corol:2}}\\
\noindent
Given~\eqref{eq:algdec}, it follows from the asymptotic properties of
$h_{(\alpha)}(t)$ and from the algebraic limit theorem
\begin{align}
\lim_{t\to \pm\infty}\frac{\dot{\Gamma}(t)}{\dot{g}_{(\alpha)}(t)}&= 
\lim_{t\to \pm\infty}\frac{\dot{\Gamma}(t)}{ t^{-m}} \frac{ t^{-m}}{\dot{g}_{(\alpha)}(t)}
=
\lim_{t\to \pm\infty}\frac{\dot{\Gamma}(t)}{t^{-m}}  
\lim_{t\to \pm \infty } \frac{ t^{-m}}{\dot{g}_{(\alpha)}(t)}
=
\lim_{t\to \pm\infty}\frac{\dot{\Gamma}(t)}{t^{-m}}  
\lim_{s\to \pm 1 } \frac{(h_{(\alpha)}(s))^{-m}}{\gamma_{(\alpha)}(s)}
\nonumber\\
&= \lim_{t\to \pm\infty}\frac{\dot{\Gamma}(t)}{t^{-m}}
\lim_{s\to \pm 1}\frac{1}{\alpha s^m}\left(
\alpha(1-s^2)^\frac{m-1}{\alpha} + 2s^2(1-s^2)^\frac{m-1-\alpha}{\alpha}
\right),\nonumber
\end{align}
 that the first transformation condition~\eqref{eq:dLdg} 
is satisfied if both $(1-s^2)$ terms above have non-negative powers, or if $0<\alpha\le m-1$.
It follows from~\eqref{eq:tc2alg} that the second transformation
condition~\eqref{eq:d2gdt2} is satisfied for $\alpha > 0$. It then follows from Proposition~\ref{l2} 
that the compactified system~\eqref{eq:dsdtalg} is $C^1$-smooth on the extended phase 
space $U\times [-1,1]$ for any $0<\alpha\le m-1$. 

The Hausdorff semi-distance $d(p(t),S^\pm)$ between a point  $p(t)=(x(t),s(t))$ on a trajectory and invariant subspace $S^\pm=\{(x,s):s=\pm 1\}$ in the extended phase space 
is the distance between $s(t)$ and $\pm 1$. Thus, algebraic decay of $d(p(t),S^\pm)$ follows from the algebraic approach of $s(t)$ 
towards $\pm1$ as $t\to\pm\infty$:
$$
d(p(t),S^\pm) = \inf_{y\in S^\pm}\,\Vert p(t) - y\Vert  = 1\mp s(t).
$$
It follows from Eq.~\eqref{eq:sLyap} that the additional Lyapunov exponent $l_{s,(\alpha)}^{\pm}$
is given by Eq.~\eqref{eq:tc2alg}.
\qed\\

Similarly, one can discuss one-sided subcases of the two-sided compactification above,
using  an example of a parametrised inverse coordinate
transformation~\eqref{eq:g+} with algebraic decay in the form
\begin{equation}
t = h_{(\alpha)}(s) = \frac{s}{(1-s)^\frac{1}{\alpha}},
\end{equation}
with the special case $\alpha=1$:
\begin{equation}
s = g(t) =\frac{t}{1 + t},\;\; t = h(s) = \frac{s}{1-s},
\end{equation}
used in the radial steady state problem~\cite{Jones1986},
or  example of parametrised inverse coordinate transformation~\eqref{eq:g-} with algebraic decay in the form
\begin{equation}
t = h_{(\alpha)}(s) = \frac{s}{(1+s)^\frac{1}{\alpha}},
\end{equation}
with the special case $\alpha=1$:
\begin{equation}
s = g(t) =\frac{t}{1 - t},\;\; t = h(s) = \frac{s}{1+s}.
\end{equation}

\subsection{Other Compactification Types and  Decay Limitation}

In Secs.~\ref{sec:ec} and~\ref{sec:ac} we constructed examples of exponential~(\ref{eq:gexp1}) and algebraic~(\ref{eq:halg1}) 
transformations, respectively, with the same type and rate of asymptotic decay 
as $t\to\pm\infty$. 
Transformations  with different  
rates or different types of asymptotic decay as $t\to+\infty$ and $t\to -\infty$ 
are also possible.
However, there is one important limitation, namely that the second
compactification condition~\eqref{eq:d2gdt2} excludes coordinate 
transformations with  faster than exponential decay; see Remark~\ref{rmk:superexp}.

In case of bi-asymptotically constant $\Gamma(t)$ whose speed $\Vert \dot{\Gamma}(t)\Vert $ has the 
same type and rate of asymptotic decay as $t\to\pm\infty$, 
the natural choice is transformation~(\ref{eq:g}) whose
derivative $\dot{g}_{(\alpha)}(t)$ has one type and rate of asymptotic decay. 
Such a transformation works also for $\dot{\Gamma}(t)$ with different rates or types of asymptotic decay.

Nonetheless, in the case of bi-asymptotically constant $\Gamma(t)$ whose speed $\Vert \dot{\Gamma}(t)\Vert $  has  
the same type but different rates of asymptotic decay as $t\to+\infty$ and $t\to -\infty$, 
one may wish to construct 
a transformation~(\ref{eq:g}) that, unlike 
examples~(\ref{eq:gexp1}) and~(\ref{eq:halg1}), has the same type but 
different rates of asymptotic decay as $t\to+\infty$ and $t\to -\infty$. 
For example, the exponential transformation~(\ref{eq:gexp1}) can be generalised to
\begin{equation}
\label{eq:gexp1gen}
t=h_{(\alpha_-,\alpha_+)}(s)= \frac{1}{\alpha_-}\ln(1 + s) - \frac{1}{\alpha_+}\ln(1-s),
\end{equation}
where the {\em compactification parameters}  $\alpha_-$ and $\alpha_+ > 0$ 
parametrise the {\em rates} of exponential decay of $\dot{g}_{(\alpha_-,\alpha_+)}(t)$ 
as $t\to-\infty$ and $t\to+\infty$, respectively.
Similarly,  algebraic transformation~(\ref{eq:halg1}) can be generalised to
\begin{equation}
\label{eq:halg1gen}
t=h_{(\alpha_-,\alpha_+)}(s)=\frac{s}{(1+s)^\frac{1}{\alpha_-}\,(1-s)^\frac{1}{\alpha_+}},
\end{equation}
where the {\em compactification parameters} $\alpha_-$ and $\alpha_+ > 0$ parametrise the {\em orders} of algebraic 
decay of $\dot{g}_{(\alpha_-,\alpha_+)}(t)$ 
as $t\to-\infty$ and $t\to+\infty$, respectively.
In the case of bi-asymptotically constant $\Gamma(t)$ whose speed $\Vert \dot{\Gamma}(t)\Vert $  has  
different types of asymptotic decay as $t\to+\infty$ and $t\to -\infty$, 
one may wish to construct a `hybrid' transformation~(\ref{eq:g}) that, unlike examples~(\ref{eq:gexp1gen}) and~(\ref{eq:halg1gen}), 
has different types of asymptotic decay as $t\to+\infty$ and $t\to -\infty$. 
For each option discussed above, the first transformation condition~\eqref{eq:dLdg} is 
satisfied by choosing the compactification parameter(s) so that  $\dot{g}_{(\alpha_-,\alpha_+)}(t)$ 
does not decay faster than $\dot{\Gamma}(t)$.

\section{Conclussion}

We have developed a general framework to facilitate stability and 
nonlinear dynamics
analysis in nonautonomous dynamical systems where the
nonautonomous terms 
decay asymptotically. 
We refer to such systems  as asymptotically autonomous dynamical systems~\cite{Markus1956}.
The theoretical work is directly motivated by a wide range of problems
from applications including pullback attractors, rate-induced 
tipping and nonlinear wave solutions, all of which fit naturally into the 
asymptotically autonomous setting. 
The analysis is based on a suitable compactification technique developed in 
Section~\ref{sec:compact}, in conjunction with a dynamical system approach to 
study invariant manifolds of saddles in the ensuing autonomous compactified 
system. 

The main obstacle to the analysis of the original nonautonomous 
system is the absence of compact invariant sets, 
such as equilibria, limit cycles or invariant tori, in its phase space.
The key idea to overcome this obstacle is twofold:
\begin{itemize}
    \item 
    { Introduce compact invariant sets to the problem:} The original nonautonomous 
    system is reformulated into an {\em autonomous compactified system} by augmenting 
    the phase space with an additional bounded but open dimension, which is then extended 
    at one or both ends by bringing in flow-invariant subspaces that carry autonomous 
    dynamics and compact invariant sets of the limit systems from infinity. 
    \item
    { Use autonomous dynamics and compact invariant sets of the limit systems 
    from infinity to analyse the original nonautonomous system:} it turns out that 
    solutions of interest form unique invariant manifolds of saddles for the limit 
    systems when embedded in the extended phase space of the compactified system.
\end{itemize}
We describe a general case with an arbitrary decay of nonautonomous 
terms and  construct the compactified system for this general case.
The main results can be summarised as follows:
\begin{itemize}
\item 
In Theorem~\ref{th2} we derive the weakest decay conditions possible 
for the existence of a $C^1$-smooth compactified system to enable the construction of invariant manifolds (stable, unstable, centre, etc.) 
in the compactified phase space.
\item
We describe the compactified system dynamics topologically in terms of  invariant manifolds of saddles in the extended phase space. We show  in Theorems~\ref{prop:saddler} and~\ref{prop:sinkr} that the manifolds of interest are unique and comprise the entire local stable  or unstable
set of the saddle, even if the saddle gains a centre direction due to non-exponential decay of the nonautonomous term(s) and the manifolds 
are not pure stable or unstable manifolds. 
\item
We construct examples of parametrised (exponential and algebraic) compactifications that can be 
implemented in actual practice.
\end{itemize}

Our general compactification framework provides an alternative tool for the analysis of asymptotically autonomous dynamical systems that complements the existing approach based on 
asymptotic equivalence of the nonautonomous system and the autonomous limit systems. 
Extending our framework to other problems such as exponential dichotomies or infinite dimensions for PDEs remains an interesting research question for future study.

\section*{Acknowledgments}
The authors would like to thank P. Ashwin, T. Carroll, F. Holland, M. Krupa and M. Rasmussen for interesting discussions.
C.X. and S.W. received funding from the the CRITICS
Innovative Training Network via the European Union’s Horizon 2020 research and innovation programme under Grant Agreement No. 643073. CJ acknowledges support from the US Office of Naval Research under grant number N00014-18-1-2204 and the US National Science Foundation under grant number  DMS-1722578.


\bibliographystyle{plain}
\bibliography{compact}

\end{document}